\documentclass[11pt,a4paper,reqno]{amsart}
\usepackage[T1]{fontenc}
\usepackage{amsmath}
\usepackage{amsthm}
\usepackage{amsfonts}
\usepackage{amssymb}
\usepackage[pdftex]{graphicx,color,epsfig}
\usepackage{amsbsy}
\usepackage{mathrsfs}
\usepackage{bbm}
\usepackage{titletoc}
\usepackage{enumerate}
\newtheorem{theorem}{Theorem}[section]

\newtheorem{proposition}[theorem]{Proposition}

\theoremstyle{remark}

\numberwithin{equation}{section}
\catcode`@=11
\def\section{\@startsection{section}{1}%
  \z@{1.5\linespacing\@plus\linespacing}{.5\linespacing}%
  {\normalfont\bfseries\large\centering}}
\catcode`@=12
\newcommand{\be}{\begin{equation}}
\newcommand{\ee}{\end{equation}}
\newcommand{\bea}{\begin{eqnarray}}
\newcommand{\eea}{\end{eqnarray}}
\newcommand{\bee}{\begin{eqnarray*}}
\newcommand{\eee}{\end{eqnarray*}}

\def\pa{\partial}

\def\RR{\mathbb{R}}

\def\ds{\displaystyle}
\def\ni{\noindent}
\def\bs{\bigskip}
\def\ms{\medskip}

\def\eps{\varepsilon}
\def\fref#1{{\rm (\ref{#1})}}
\def\pref#1{{\rm \ref{#1}}}

\def\calE{{\mathcal E}}
\def\calO{{\mathcal O}}

\catcode`@=11
\def\supess{\mathop{\operator@font Sup\,ess}}
\catcode`@=12
\def\un{{\mathbbmss{1}}}

\def\RR{\mathbb{R}}

\def\lim{\mathop{\rm lim}}

\def\un{{\bf 1}}

\def\fref#1{{\rm (\ref{#1})}}
\def\pref#1{{\rm \ref{#1}}}
\def\ds{\displaystyle}
\def\ni{\noindent}
\def\bs{\bigskip}
\def\ms{\medskip}

\def\pa{\partial}

\newcommand{\dt}{\partial_t}

\newcommand{\grdx}{\nabla_x}

\newcommand{\cint}[1]{\left\langle #1 \right\rangle}

\def\fref#1{{\rm (\ref{#1})}}

\title[AP schemes including boundary layers]
{Micro-macro schemes for kinetic equations including boundary layers}
\author[M. Lemou]{Mohammed Lemou}
\address{CNRS and IRMAR, Universit\'e de Rennes 1, France}
\email{mohammed.lemou@univ-rennes1.fr}
\author[F. M\'ehats]{Florian M\'ehats}
\address{IRMAR, Universit\'e de Rennes 1, France}
\email{florian.mehats@univ-rennes1.fr}
\begin{document}


\maketitle
\begin{abstract}
We introduce a new micro-macro decomposition of collisional  kinetic equations in the specific case of the diffusion limit, which naturally incorporates the incoming boundary conditions. The idea is to write the distribution function $f$ in all its domain as the sum of an equilibrium adapted to the boundary (which is not the usual equilibrium associated with $f$) and a remaining kinetic part. This equilibrium is defined such that its incoming velocity moments coincide with the incoming velocity moments of the distribution function. A consequence of this strategy is that no artificial boundary condition is needed in the micro-macro models and the exact boundary condition on $f$ is naturally transposed to the macro part of the model. This method provides an 'Asymptotic preserving' numerical scheme which generates a very good approximation of the space boundary values at the diffusive limit, without any mesh refinement in the boundary layers. Our numerical results are in very good agreement with the exact so-called Chandrasekhar value, which is explicitely known in some simple cases. 
\end{abstract}

\titlecontents{section} 
[1.5em] 
{\vspace*{0.1em}\bf} 
{\contentslabel{2.3em}} 
{\hspace*{-2.3em}} 
{\titlerule*[0.5pc]{\rm.}\contentspage\vspace*{0.1em}}

\vspace*{-2mm}

\section{Introduction} \label{intro}

The development of numerical methods to solve multiscale kinetic equations has been the subject of active research in the past years, with applications in various fields: plasma physics, rarefied gas dynamics, aerospace engineering, semiconductors, radiative transfert, ... 
The general problem is to construct numerical schemes that are able to capture the properties of the various scales in the considered system, while the numerical parameters remain as independent as possible of the stiffness character of these scales. 

There is a huge literature dealing with the construction of such schemes, and in particular it is now a usual challenge to design the so-called {\em Asymptotic Preserving} (AP) \cite{jin} numerical schemes for kinetic equations which are consistent with the kinetic model for all positive value of the Knudsen number $\eps$, and degenerate into consistent schemes with the asymptotic models (compressible Euler and Navier-Stokes, diffusion, etc) when $\eps \to 0$. Among a long list of works on this topic, we can mention \cite{larsen-morel,larsen-morel-miller,jin-levermore1,jin-levermore2,GJL,jin-tang,guermond-kanschat} for stationary problems and \cite{JPT,jin-pareschi,Klar1,Klar2,Klar3,klar4,gosse-toscani,LM-AP1,liu-mieussens,jin-filbet,goudon,pareschi-recent,lafitte} for time dependent problems.

Here, we are interested in boundary value problems and our aim is to develop a strategy to construct numerical schemes that are AP in the usual sense and are also able to deal with the space boundary conditions and kinetic boundary layers. In this paper, we focus on the diffusive scaling but our strategy could be extended to other asymptotics. We emphasize that in general, the known numerical schemes induce at the limit $\eps\to 0$ some numerical boundary conditions which are different from the theoretical limiting boundary values. In most cases, this does not affect only the boundary layer, but also the solution inside the domain.

Boundary conditions in the construction of AP schemes for transient kinetic equation in the diffusive scaling have been taken into account in several works. For instance, in \cite{Klar1,Klar3,klar4}, Klar uses an approximation of the solution to the Milne problem which is injected as a boundary condition for the evolution problem. The boundary value used by Jin, Pareschi and Toscani in \cite{JPT} is also derived from the asymptotic analysis of the problem as the Knudsen number $\eps$ goes to zero. Of course, taking the boundary condition from the limiting model cannot provide a good approximation of the problem for all regime, in particular for large values of $\eps$. Our goal here is to develop a new approach in which the boundary conditions are not extracted from the limiting model but are directly derived from the original kinetic boundary condition, at any time. Our numerical results will be compared with those obtained with the approaches in \cite{Klar1,JPT,LM-AP1}.

To this purpose, our starting idea is to extend the micro-macro decomposition of \cite{BLM,LM-AP1,Lem-note1} in order to incorporate in a natural way the exact inflow boundary conditions. The general micro-macro decomposition method consists in writing the kinetic equation as a system coupling a macroscopic part (say a Maxwellian) and a remaining kinetic part.  The main advantage of this method is its robustness and easy adaptability, since it can be applied to various asymptotics (fluid, diffusion, high-field, etc) and to a large class of collision operators.  However, the general problem of space boundary condition and boundary layers has not been completely solved in this approach.  When incoming boundary conditions on the distribution function  are prescribed, it is clear that this cannot be translated into a boundary condition on the macro part (say the total density) in an exact way since the distribution function is only known for incoming velocities. Such boundary conditions for both macro and micro parts of the distribution function are needed for the computation of fluxes. In this case, an artificial boundary condition of Neumann type is used in \cite{LM-AP1} and leads to some unsatisfactory approximation of the Chandrasekhar value at the diffusive limit. 

Our first aim is to develop a new micro-macro decomposition of collisional kinetic equations which naturally incorporates the exact space boundary conditions (BC). The second task in this work is  to show how to use this new formulation at the boundary in order to capture the right boundary conditions and boundary layers in the limit of  a small Knudsen number.

To achieve this goal we proceed as follows. We consider a kinetic equation with a linear collision operator in a diffusive scaling on a bounded domain, subject to inflow boundary conditions. Our first idea is to decompose the distribution function $f$ in its domain as the sum of a Maxwellian part adapted to the boundary (which is not the usual Maxwellian associated with $f$) and a remaining kinetic part. This Maxwellian is defined on the whole domain such that its 'incoming' velocity moments coincide with the 'incoming' velocity moments of the distribution function. In this way, the exact boundary condition can be directly injected into the micro-macro system. The second idea is to use a suitable approximation of the space derivative of the density at the boundary which garantees a balance between the incoming and outgoing fluxes at this boundary. Important consequences of this strategy are the following. No artificial boundary condition is needed in the micro/macro models and the exact boundary conditions on $f$ are naturally shared by the macro part and the kinetic part.  A further fundamental property of the strategy developed here is that it provides AP schemes inside the physical domain and a very good approximation in the space boundary layers. Indeed, at the end of this paper, we provide various numerical tests which illustrate this property. We emphasize that the numerical boundary value given by our scheme when $\eps$ is small is very close to the theoretical boundary value derived from the Chandrasekhar function, without injecting this value in our scheme. The strategy developed in this paper has been announced in a preliminary note \cite{note}, where we showed that this method can also be applied in principle for a hydrodynamic scaling. We will explore in the future several extensions of this work, in particular the possibility of dealing with more complex collision operators. Our approach shares some similarities with the method developed in \cite{dubroca,dubroca2} where a partial moment system with entropy closure was used to approximate the original kinetic equation.

The paper is organized as follows. In Section \ref{sect2}, the general strategy is described at the continuous level. We introduce a new micro-macro decomposition which is adapted to boundary value problems for kinetic equations. In Section \ref{sectionscheme}, we use this decomposition to design a numerical scheme which is AP in the diffusion limit. In a first step (Subsection \ref{sectattempt}), we show how this formulation allows to inject directly the inflow boundary condition in the numerical method, avoiding in this way any artificial boundary condition. In a second step (Subsection \ref{numscheme}), we introduce a specific (still consistent) discretization of the spatial derivative of the density in order to correctly capture the boundary condition and the boundary layer. The properties of our scheme are summarized in Proposition \ref{propschema2}. In Section \ref{sectionnum}, we provide various numerical tests to validate our approach. We test different boundary conditions (with or without boundary layers) and different linear collision operator (with local and non local collision kernels). Finally, we end the paper by a conclusion and some perspectives.





\section{A boundary matching micro-macro decomposition}
\label{sect2}

\subsection{The diffusion limit of kinetic equations in bounded domains}

Let $\Omega$ be a domain in the position space $\RR^d$ and let $V$ be a domain in the velocity space $\RR^d$, $V$ being endowed with a measure $d\mu$. At any point $x\in\pa\Omega$ of the boundary of $\Omega$, we denote by $n(x)$ the unitary outgoing normal vector to $\pa\Omega$. We then consider the following linear transport equation written in a diffusive scaling:
\begin{equation}  \label{eq-transp}
\eps \pa_t f + v\cdot \nabla_x f = \frac{1}{\eps}L f ,   \qquad t > 0,  \qquad
  (x,v) \in \Omega \times V,
\end{equation}
with initial condition
\be
\label{initialdata}
\qquad f(0,x,v)=f_{init}(x,v),  \qquad (x,v) \in \Omega \times V
\ee
and with inflow boundary conditions
\be
\label{BC-f}
f(t,x,v)= f_{b}(t,x,v), \qquad t>0, \quad (x,v) \in \partial \Omega \times V\ \mbox{such that} \ v\cdot n(x) < 0.
\ee
The unknown $f$ is the distribution function of the particles that depends on time $t>0$, on position $x\in\Omega $, and on velocity
$v\in V$.  

In \fref{eq-transp}, the linear collision operator $L$ acts on the velocity dependence of $f$ and describes the interactions of particles with the medium. This operator relaxes the system of particles to an equilibrium $\mathcal E(v)$, which is supposed to be a positive and even function.
For all distribution function $h$, we shall denote
\be
\label{crochetV}
\cint{h}_V=\frac{\int_Vh(v)dv}{\int_V\mathcal E(v)dv}.
\ee
Note in particular that $\cint{\mathcal E}_V=1$ and $\cint{v\mathcal E}_V=0$. We assume that the collision operator $L$ is non-positive and self-adjoint in $L^2(V,\mathcal E^{-1}\, d\mu)$, with null space and range given by
 $${\mathcal N}(L)=\mbox{Span}\lbrace\mathcal E\rbrace$$ 
and
$${\mathcal R}(L) = ({\mathcal N}(L))^\bot = \lbrace f \text{ such that } \cint{f}_V=0 \rbrace.$$
Hence, an important property of $L$ is the fact that it is invertible on $\mathcal R(L)$, we shall denote its pseudo-inverse by $L^{-1}$. Additionally, we assume the invariance of $L$ under orthogonal transformations of $\RR^d$. The parameter $\eps>0$ measures a dimensionless mean free path of the particles, or also the inverse of the dimensionless observation time.

When $\eps$ goes to 0 in \fref{eq-transp}, $f$ goes formally to an equilibrium state $f_0(t,x,v)=\rho_0(t,x) \mathcal E(v)$. The diffusion limit is the equation satisfied by the limiting density $\rho_0$. This equation is a diffusion equation and has been obtained in various situations thanks to asymptotic expansion in $\eps$ (Hilbert or Chapman-Enskog expansion), see e.g. \cite{larsen-keller,BLP,BSS,degond-masgallic,poupaud}. We have
\begin{equation}  \label{eq-diff}
\dt \rho_0 - \nabla_x \cdot (\kappa \nabla_x \rho_0) = 0 \qquad  \mbox{with}\qquad \kappa=-\cint{vL^{-1}(v\mathcal E)}_V>0,
\end{equation}
with a Dirichlet boundary condition
\be
\label{bc-rho}
\rho_0(t,x)=\rho_b(t,x),\qquad t>0,\quad x\in \pa\Omega.
\ee
The limit value $\rho_b$ is usually obtained by a boundary layer analysis and its computation involves the resolution of a kinetic half-space problem, the Milne problem associated with the collision operator $L$. For $x\in \pa\Omega$, let $\chi^{t,x}(y,v)$ be the bounded solution of the half-space problem
\bea
\label{milne}
-v\cdot n(x)\pa_y\chi^{t,x}&=&L\chi^{t,x},\qquad \qquad y>0,\quad v\in V,\nonumber\\
\chi^{t,x}(0,v)&=&f_b(t,x,v),\qquad v\cdot n(x)<0,
\eea
in which $t$ and $x$ are parameters. Then, under reasonnable conditions of $f_b$, this problem is well-posed \cite{BLP,BSS,golse,poupaud} and we have
\be
\label{exact}
\lim_{y\to +\infty}\chi^{t,x}(y,v)=\rho_b(t,x)\calE(v).
\ee
Notice that, in the special case where the incoming data is an equilibrium state, i.e. when $f_b(t,x,v)=\rho_b(t,x)\calE(v)$, one has $\chi^{t,x}(y,v)=\rho_b(t,x)\calE(v)$ and is independent of $y$. In this case, there is no boundary layer near $\pa\Omega$ as $\eps\to 0$. For general incoming distribution functions, the question of determining $\rho_b$ requires the resolution of the half-space problem \fref{milne}, which may be computationaly demanding. Our aim here is to construct a numerical method which gives a good approximation of the boundary layer without solving the half-space problem \fref{milne}.

\bs
In order to check the efficiency of our approach, a validation test will be the following. We will compare our numerical results with the exact Dirichlet value $\rho_b$ defined by \fref{exact} in the diffusion limit $\eps\to 0$, in special cases where this value is well-known. The most simple case is the one-group transport equation in slab geometry, which is such that $d=1$, $\Omega=[0,1]$, the velocity set is $V=[-1,1]$, $d\mu=\frac12 dv$ and $Lf=\langle f\rangle_V -f$. Then $\rho_b$ is given by
\be
\label{chandra1}
\rho_b(t,0)=\int_0^1K_0(v)f_b(t,0,v)dv,\qquad \rho_b(t,1)=\int_{-1}^0K_0(-v)f_b(t,1,v)dv,
\ee
where 
\be
\label{chandra}
K_0(v)=\frac{\sqrt{3}}{2}vH(v),
\ee
the function $H$ being the so-called Chandrasekhar function \cite{chandrasekhar,BSS,dautray-lions}, see also \cite{case-zweifel}, which can be computed numerically from the following formula:
\be
\label{pointfixe}
H(v)=1+v\frac{H(v)}{2}\int_0^1\frac{H(w)}{v+w}\,dw.
\ee
A good approximation of the function $K_0(x)$ is usually given by 
\be
\label{polynome1}
K_0(v)\simeq K_1(v)=\frac{3}{2}v^2+v,
\ee
see \cite{larsen-morel} for instance.

\subsection{A new micro-macro decomposition}

Let us first recall the micro-macro decomposition introduced in \cite{LM-AP1} and \cite{Lem-note1}.  We decompose
\be
\label{mm1}
f= \rho \mathcal E + g_1,
\ee
with $$\rho(t,x)= \cint{f}_V\quad \mbox{and}\quad g_1=f-\rho \mathcal E$$ ($g_1$ is not necessarily small). We denote by $\Pi$ the orthogonal projector in $L^2(\mathcal E^{-1}d\mu)$ onto the nullspace of $L$, i.e. $$\Pi \phi= \cint{\phi}_V  \mathcal E.$$ Inserting the decomposition \fref{mm1} into the kinetic equation \fref{eq-transp} and applying  $\Pi$ and $I -\Pi$  successively ($I$ being the identity operator), one gets
\bee
&&\dt \rho + \frac{1}{\eps}  \grdx \cdot
  \cint{ v g_1 }_V =  0, \\  
&&   \dt g_1 + \frac{1}{\eps} (I - \Pi)(v \cdot \grdx g_1) = \frac{1}{\eps^2}
    \big[ L g_1 - \eps  \mathcal Ev\cdot \grdx \rho) \big].
\eee
We now emphasize that in this formulation, the space boundary condition on $\rho$ and $g_1$ are not known because they cannot be inferred from the boundary conditions on $f$ in general.  Indeed, \fref{BC-f} only gives the function $f$ at the boundary for incoming velocities, i.e. such that $v\cdot n(x) < 0$, and it is therefore clear that the values of $\rho(t,x) =\cint{f}_{V}$ cannot be determined a priori on the boundary $\partial \Omega$ with the only knowledge of $f_b$. Consequently, the micro-macro decomposition \fref{mm1} appears as a method for capturing the diffusion limit inside the domain only. The numerical scheme constructed in \cite{LM-AP1} is Asymptotic Preserving in a strict sense (i.e. including the boundary) only for incoming distribution data at the equilibrium $f_b(t,x,v)=\rho_b(t,x)\calE(v)$, when there is no boundary layer near $\pa \Omega$. Indeed, only a rough approximation of the boundary layer is obtained for non equilibrium boundary conditions.

Our aim in this work is to develop a new micro-macro decomposition  which is able to incorporate the boundary conditions in an exact way. This will allow us in the next section to construct a numerical method with natural boundary conditions, that approximate very well the diffusion limit even for incoming data that are not an equilibrium state.

To this purpose, let us introduce a few further notations. We consider a function $\omega(x,v)$ which extends the function $n(x)\cdot v$ (defined on $\pa\Omega\times V$) to the whole domain $\Omega\times V$. Let us give explicit examples of $\omega(x,v)$ for specific geometries. For the unit ball centered at the origin, one can take $\omega(x,v)=x\cdot v$. For a half plane $x_{1}>0$, one can choose $\omega(x,v)=(-v_1, 0, ...,0)$.  In dimension one, for the interval $[0,1]$, one can take $\omega(x,v)=(2x-1)v$. Then, for all $x\in \Omega$, we split the velocity space $V$ into two parts according to the sign of this function $\omega(x,v)$: 
\be
\label{VOpm}
\begin{array}{ll}
V_{-}(x)= \{v\in V, \ \omega(x,v)<0 \},  &\qquad \    \ V_{+}(x)= V\backslash V_{-}(x).
\end{array}
\ee
For all function $h(v)$, we will denote
$$ \cint{h}_{V_{-}}=\frac{\int_{V_{-}} h(v)  d\mu}{\int_{V_{-}} \calE(v) d\mu},\qquad \cint{h}_{V_{+}}=\frac{\int_{V_{+}} h(v)  d\mu}{\int_{V_{+}} \calE(v) d\mu}
$$
and
$$\Pi_{V_-} h=\cint{h}_{V_{-}}\calE.$$
The idea is now the following. Instead of looking for an equation on $\rho$ as usual, we seek  an equation on the following 'boundary matching' density:
\be 
\label{rhobar-def}\overline \rho(t,x) =\cint{f(t,x,\cdot)}_{V_{-}}
\ee
 and perform the corresponding micro-macro decomposition:
\be
\label{mm2}
  f= \overline \rho \mathcal E + g=\Pi_{V_-} f+g.
 \ee
In particular, we have $\cint{g}_{V_-}=0$. When $\eps \to 0$, we know that the solution $f$ of \fref{eq-transp} is, at least formally, close to $\rho \mathcal E$ except in initial or boundary layers. Therefore, since we have $\cint{\rho \mathcal E}_{V_{-}}=\rho,$ the moment $\overline \rho$ will also be close to $\rho$ for small $\eps$.  This shows that the new decomposition \fref{mm2} is still a decomposition of $f$ into an asymptotic part (macro part) and a kinetic part (micro part).
In order to derive the system of equations satisfied by $\overline \rho$ and $g$, we first integrate \fref{eq-transp} on $V_{-}$ and get
\be
\label{equ-rhobar}
\partial_{t}\overline \rho + \frac{1}{\eps}\cint{v\mathcal E}_{V_{-}}\cdot  \nabla _{x}\overline\rho + \frac{1}{\eps}\cint{v\cdot \nabla _{x}g}_{V_{-}}= \frac{1}{\eps^2}\cint{Lg}_{V_{-}}.
\ee
Substracting \fref{equ-rhobar} from equation \fref{eq-transp} on $f$, we obtain the equation on $g$:
\be
\label{equ-g1}
\partial_{t}g + \frac{1}{\eps}(I-\Pi_{V_-} )( v\cdot \nabla_{x}g)+ \frac{1}{\eps}(I-\Pi_{V_-} )(v \calE)\cdot \nabla_{x}\overline \rho  = \frac{1}{\eps^2}(I-\Pi_{V_-})(Lg).
\ee
Finally, we remark that the system \fref{equ-rhobar}, \fref{equ-g1} can be replaced by a more convenient (and still equivalent) system in terms of $\rho =\cint{f}_V$ and $\ds g= f- \overline \rho \mathcal E$ as follows:
\bea
\label{MM-rho}
&&\hspace*{-1.1cm}\partial_{t} \rho +  \frac{1}{\eps}\cint{v\cdot \nabla _{x}g}_V=0,\\
\label{MM-g}
&&\hspace*{-1.1cm}\partial_{t}g + \frac{1}{\eps}(I-\Pi_{V_-} )( v\cdot \nabla_{x}g)+ \frac{1}{\eps}(I-\Pi_{V_-} )(v \calE)\cdot \nabla_{x}\overline \rho  = \frac{1}{\eps^2}(I-\Pi_{V_-})(Lg).
\eea
Note that $\overline \rho$ is linked to $\rho$ and $g$ by the relations 
\be
\label{rhobarrho} \overline \rho  = \rho - \cint{g}_V, \qquad f= \rho \calE+g -\cint{g}_V\calE.
\ee
One main interest of this new  micro-macro formulation  \fref{mm2} is the following: the fluxes involved in \fref{MM-rho} and \fref{MM-g} only concern the quantities $g$ or $\overline \rho$, and not $\rho$. This means that numerical schemes of this formulation would only need the values of $\overline \rho$ and $g$ at the space boundary which are completely known from the original boundary condition  \fref{BC-f} on $f$:
\bea
\label{mm-bc1}
&&\hspace*{-1.2cm}\overline{\rho}(t,x,v)=\cint{f_b(t,x,\cdot)}_{V_{-}},\qquad \qquad \qquad \quad \qquad  t>0, \quad (x,v) \in \partial \Omega,\\
\label{mm-bc2}
&&\hspace*{-1.2cm}g(t,x,v)= f_{b}(t,x,v)-\cint{f_b(t,x,\cdot)}_{V_{-}}\calE(v), \qquad t>0, \quad (x,v) \in \partial \Omega \times V_-.
\eea
Note that a similar decomposition as \fref{mm2} can be performed for a fluid scaling as well, see \cite{note}.

\section{The numerical method}
\label{sectionscheme}

In this section, we construct in dimension $d=1$ a numerical scheme for \fref{MM-rho}, \fref{MM-g}, see Proposition \ref{propschema2}, which is uniformly stable in the limit $\eps\to 0$ and provides a good approximation of the boundary layers, without introducing any artificial boundary condition. To this purpose, we proceed in two steps. First, we present in subsection \ref{sectattempt} a preliminary numerical scheme, see Proposition \ref{propschema1}, which has the desired uniform stability and does not use any artificial boundary condition. However, it turns out that this scheme does not give a good approximation of the limiting boundary condition. To remedy this problem, we then introduce our numerical scheme in subsection \ref{numscheme} where a specific discretization is done near the boundary. We emphasize that this scheme is able to reproduce an accurate approximation of the boundary layer without any mesh refinement.

\bs
Let us introduce a few notations. The domain in space is the interval $[0,1]$. The incoming boundary conditions are
\bea
\label{bc1D0}
f(t,0,v)=f_\ell(v),\mbox{ for }v\in V_-(0)=\{v\in V,\,v>0\},\\
\label{bc1D1}
\qquad f(t,1,v)=f_r(v),\mbox{ for }v\in V_-(1)=\{v\in V,\,v<0\},
\eea
the data $f_\ell$ and $f_r$ being independent of time here for simplicity. The function $\omega(x,v)$ is a smooth function satisfying
\be
\label{ome}\omega(0,v)=-v,\qquad \omega(1,v)=v.
\ee
Let $\Delta t$ be a constant time step and $\Delta x=\frac{1}{N+1}$ be a uniform space step. We set $t_n=n\Delta t$ and consider staggered grids $x_i=i\Delta x$ for $i=0,\ldots,N+1$ and $x_{i+1/2}=(i+1/2)\Delta x$ for $i=0,\ldots,N$. We will construct approximations of the densities $\rho^n_i\simeq \rho(t_n,x_i)$ and $\overline\rho^n_i\simeq \overline\rho(t_n,x_i)$ on the grid $\{x_i\}_{0\leq i\leq N+1}$ and approximations of the remainder $g^n_{i+1/2}(v)\simeq g(t_n,x_{i+1/2},v)$ on the grid $\{x_{i+1/2}\}_{0\leq i\leq N}$. The velocity variable $v$ is kept continuous in this section, and its discretization will be made precise later on. At the time $t=0$, the initial condition \fref{initialdata} induces naturally the initialization
\be
\label{initialdisc}
\rho^0_i=\cint{f_{init}(x_i,\cdot)}_V,\quad \overline\rho^0_i=\cint{f_{init}(x_i,\cdot)}_{V_-},\quad g^0_{i+1/2}=(I-\Pi_{V^-})f_{init}(x_{i+1/2},\cdot).
\ee
Note that, by construction, the function $\overline \rho$ is known at the boundary $x=0$, $x=1$ and can be deduced from the boundary conditions \fref{bc1D0}, \fref{bc1D1}:
\be
\label{incbar}\overline \rho_0=\frac{\int_{V_-(0)} f_\ell(v)  d\mu}{\int_{V_-(0)} \calE(v) d\mu},\qquad \overline \rho_{N+1}=\frac{\int_{V_-(1)} f_r(v)  d\mu}{\int_{V_-(1)} \calE(v) d\mu}.
\ee
At the boundary, the function $g$ is also known for entering velocities. We introduce two additional unknowns $g_0^{n}(v)$ and $g_{N+1}^{n}(v)$ for $v\in V$. The incoming boundary conditions yield
\bea
\nonumber
&&g^n_0(v)=f_\ell(v)-\overline \rho_0\calE\qquad \quad \mbox{ for }v\in V_-(0)\\
\label{inc}
&&g^n_{N+1}(v)=f_r(v)-\overline \rho_{N+1}\calE\quad \mbox{ for }v\in V_-(1).
\eea
For outgoing velocities, these functions will be provided by our numerical scheme, when their values are required, and no artificial condition will be needed. This is explained in Subsection \ref{numscheme}.

\subsection{A first attempt without specific treatment at the boundary}
\label{sectattempt}
We assume that we know $\rho^n_i$ and  $\overline\rho^n_i$ for $i=1,\ldots,N$, and $g^n_{i+1/2}$ for $i=0,\ldots N$. We proceed by recursion and explain how to compute these same quantities at time $t_{n+1}$. First, for $i=0,\cdots,N,$ we discretize \fref{MM-g} by a semi-implicit upwind scheme:
\bea
\hspace*{-1cm}\frac{g_{i+1/2}^{n+1}-g_{i+1/2}^n}{\Delta t}&+&\frac{1}{\eps}(I-\Pi_{V_-})\left(v^+\frac{g_{i+1/2}^n-g_{i-1/2}^n}{\Delta x}+v^-\frac{g_{i+3/2}^n-g_{i+1/2}^n}{\Delta x}\right)\nonumber\\
\hspace*{-1cm}&+&\frac{1}{\eps}(I-\Pi_{V_-})(v\mathcal E)\left(\frac{\overline \rho^n_{i+1}-\overline \rho^n_{i}}{\Delta x}\right)=\frac{1}{\eps^2}(I-\Pi_{V_-})(Lg_{i+1/2}^{n+1}),
\label{schemag1}\eea
where $v^+=\max(v,0)$ and $v^-=\min(v,0)$. Note that the source term involving $\pa_x \overline \rho$ is discretized by a centered finite difference scheme.
We observe that $g_{-1/2}^n$ and $g_{N+3/2}^n$ for incoming velocities can be extracted from the relations
\be
\label{schemainter1}
g_0^n=\frac{g_{-1/2}^n+g_{1/2}^n}{2},\qquad g_{N+1}^n=\frac{g_{N+1/2}^n+g_{N+3/2}^n}{2},
\ee
since $g_{1/2}^n$ and $g_{N+1/2}^n$ are known and since $g_0^n$ and $g_{N+1}^n$ are known for incoming velocities by \fref{inc}.
Note also that $\overline \rho_0$ and $\overline \rho_{N+1}$ are known by \fref{incbar}.
Therefore, the scheme \fref{schemag1} enables to compute all the desired values $g^{n+1}_{i+1/2}$ for $i=0,\ldots N$ (see below in the proof of Proposition \ref{propschema1}).

It remains to compute $\rho^n_i$ and  $\overline\rho^n_i$ for $i=1,\ldots,N$. To that purpose, we consider the following discretization of \fref{MM-rho}. For $i=1,\ldots,N$, we set
\be
\label{schemarho1}
\frac{\rho^{n+1}_i-\rho^n_i}{\Delta t}+\frac{1}{\eps}\cint{v\frac{g_{i+1/2}^{n+1}-g_{i-1/2}^{n+1}}{\Delta x}}_V=0.
\ee
These relations clearly provide all the desired values $\rho^n_i$ for $i=1,\ldots,N$. Then we compute all the values $\overline \rho^{n+1}_i$ according to \fref{rhobarrho}, by setting
\be
\label{schemarhobar1}\overline \rho^{n+1}_i= \rho^{n+1}_i-\cint{\frac{g_{i+1/2}^{n+1}+g_{i-1/2}^{n+1}}{2}}_V
\ee
for $i=1,\ldots,N$. The recursion is now complete and can be iterated. Observe that, in this first attempt, we do not use the outgoing values of $g_0^n$ and $g_{N+1}^n$ in order to compute $\rho^{n+1}_i$ and $g^{n+1}_{i+1/2}$ inside the domain.

We summarize in the following proposition some properties of this scheme.
\begin{proposition}[Formal]
\label{propschema1}
The numerical scheme \fref{schemag1}, \fref{schemainter1}, \fref{schemarho1}, \fref{schemarhobar1} with initial and boundary conditions \fref{initialdisc}, \fref{incbar}, \fref{inc} determines all the values $\rho^n_i$ and  $\overline\rho^n_i$ for $i=1,\ldots,N$, and $g^n_{i+1/2}$ for $i=0,\ldots N$, at any time interation $n$. This scheme is stable under the CFL condition
$$\Delta t\leq C\left(\Delta x^2+\eps\Delta x\right),$$
where $C$ is a constant independent of $\eps$, $\Delta t$ and $\Delta x$. For all fixed $\eps>0$, this scheme is consistent with the initial and boundary value problem \fref{eq-transp}, \fref{initialdata}, \fref{BC-f}. Moreover, for all fixed $\Delta t$ and $\Delta x$, the scheme degenerates as $\eps \to 0$ into a discretization of the diffusion equation \fref{eq-diff}. More precisely, for $i=1,\ldots, N$, we have  $\rho_i^n\to  \widetilde \rho_i^n$ as $\eps \to 0$, and $\widetilde \rho^n_i$ is defined by the scheme
$$
\frac{\widetilde \rho^{n+1}_i-\widetilde \rho^n_i}{\Delta t}-\kappa\frac{\widetilde \rho^n_{i+1}-2\widetilde \rho^n_i+\widetilde \rho^n_{i-1}}{\Delta x^2}=0,\quad \mbox{for}\quad i=1,\ldots,N,
$$
where $\kappa>0$ is given in \fref{eq-diff} and where the boundary values are
$$\widetilde \rho_0^n=2\frac{\langle v L^{-1}(I-\Pi)(v^+f_\ell)\rangle_V}{\langle v L^{-1}(v\calE)\rangle_V},\qquad \widetilde \rho_{N+1}^n=2\frac{\langle v L^{-1}(I-\Pi)(v^-f_r)\rangle_V}{\langle v L^{-1}(v\calE)\rangle_V},$$
for all $n$. We recall that $\Pi h=\langle h\rangle_V \calE$.
\end{proposition}
In the particular case where $V=[-1,1]$, $\calE=1$, $d\mu=\frac{1}{2}dv$ and $L=\Pi-I$, we obtain $\kappa=\frac{1}{3}$ and 
$$\widetilde \rho_0=\int_0^1K_2(v)f_\ell(v)dv,\qquad \widetilde \rho_{N+1}=\int_{-1}^0K_2(-v)f_r(v)dv.$$
with
\be
\label{polynome2}
K_2(v)=3v^2.
\ee
For instance, if $f_\ell(v)=v$, we get the value $\widetilde \rho_0=0.75$ as boundary value for the limiting diffusion equation, whereas the right value can be computed from \fref{chandra1} and \fref{chandra}: $\rho_b(0)=\frac{\sqrt{3}}{2}\int_0^1 v^2H(v)dv\simeq 0.7104$. The value obtained by this scheme should clearly be improved. This is the goal of the next section. Note that the same boundary value 0.75 is obtained by the schemes proposed in \cite{LM-AP1} and \cite{JPT}, which both involve some artificial boundary conditions.
We emphasize that, by construction, our scheme only involves the natural boundary conditions \fref{incbar} and \fref{inc} derived from $f_\ell$ and $f_r$ and does not need any artificial condition.

\bs
\ni
{\em Proof of Proposition \pref{propschema1}.}
\nopagebreak

\ms
\ni
We first check that all the values are uniquely defined by the scheme. The algorithm has been described above and it only remains to show that  \fref{schemag1} allows to compute $g^{n+1}_{i+1/2}$ in terms of quantities at time $n$. To this aim, we simply remark that the calculation of $g^{n+1}_{i+1/2}$ is done by solving a linear system of the form
\be
\label{inversion}
\left(\frac{\eps^2}{\Delta t} I-\widetilde L\right)g=h,\quad \mbox{with}\quad \widetilde L=(I-\Pi_{V_-})L(I-\Pi_{V_-}).
\ee
Since the operator $-L$ is self-adjoint and nonnegative, so is the operator $-\widetilde L$ and then the operator $\frac{\eps^2}{\Delta t} I-\widetilde L$ is invertible, for $\eps>0$.

For the stability CFL condition, a similar proof as in \cite{LM-AP1, liu-mieussens} can be done. We only sketch the main arguments of this proof and refer to these works for the details. Roughly, when $\eps=\calO(1)$, the convection term in the equation \fref{schemag1} prevails and induces the standard stability condition $\Delta t\leq C\eps \Delta x$. When $\eps \ll 1$, the system behaves as the diffusion limiting equation (see below) and  this results in the usual stability condition $\Delta t\leq C\Delta x^2$. In the general case, the CFL condition is expressed as a combination of these two limiting situations.

Now, we analyze formally the scheme when $\eps\to 0$, the parameters $\Delta t$, $\Delta x$ being fixed.  First, we remark that the operator $(I-\Pi_{V_-})L$ is invertible from $\{f:\,\langle f\rangle_{V_-}=0\}$ to itself. More precisely, one has 
$$(I-\Pi_{V_-})Lg=h, \quad \mbox{with}\quad g\in \{f:\,\langle f\rangle_{V_-}=0\},$$
if, and only if,
$$g=(I-\Pi_{V_-})L^{-1}(I-\Pi)h,\quad \mbox{with}\quad h\in \{f:\,\langle f\rangle_{V_-}=0\}.$$
Now, by induction, it is clear from \fref{schemag1} that $g_{i+1/2}^n=\calO(\eps)$, for $i=0,\ldots,N$ and $n\geq 1$, and then, from \fref{schemarho1}, we deduce that $\rho_i^n=\calO(1)$ and $\overline \rho_i^n=\calO(1)$ for $i=1,\ldots,N$ and $n\geq 1$. Reinserting these estimates in \fref{schemag1}, we get
\be
\label{g32}
g_{i+1/2}^{n+1}=\eps (I-\Pi_{V_-})L^{-1}(v\calE)\left(\frac{\overline \rho^n_{i+1}-\overline \rho^n_{i}}{\Delta x}\right)+\calO(\eps^2)
\ee
for $i=1,\ldots,N-1$, and
$$g_{1/2}^{n+1}=\eps (I-\Pi_{V_-})L^{-1}\left[v\calE\left(\frac{\overline \rho^n_{1}-\overline \rho_{0}}{\Delta x}\right)-\frac{2}{\Delta x}(I-\Pi)\left(v^+g_0^n\right)\right]+\calO(\eps^2)$$
$$g_{N+1/2}^{n+1}=\eps (I-\Pi_{V_-})L^{-1}\left[v\calE\left(\frac{\overline \rho_{N+1}-\overline \rho^n_{N}}{\Delta x}\right)+\frac{2}{\Delta x}(I-\Pi)\left(v^-g_{N+1}^n\right)\right]+\calO(\eps^2).$$
Inserting these equations into \fref{schemarho1} gives, after some computations,
$$
\frac{\rho^{n+1}_i-\rho^n_i}{\Delta t}-\kappa\frac{\rho^n_{i+1}-2\rho^n_i+\rho^n_{i-1}}{\Delta x^2}=\calO(\eps),\quad \mbox{for}\quad i=2,\ldots,N-1
$$
and
$$
\frac{\rho^{n+1}_1-\rho^n_1}{\Delta t}-\kappa\frac{\rho^n_{2}-2\rho^n_1+\tilde \rho^n_{0}}{\Delta x^2}=\calO(\eps),
$$
$$
\frac{\rho^{n+1}_N-\rho^n_N}{\Delta t}-\kappa\frac{\tilde \rho^n_{N+1}-2\rho^n_{N}+\rho^n_{N-1}}{\Delta x^2}=\calO(\eps),
$$
the quantities $\tilde \rho^n_{0}$ and $\tilde \rho^n_{N+1}$ being defined in Proposition  \ref{propschema1}.
\qed

\subsection{The numerical scheme}
\label{numscheme}
In this section, we construct a numerical scheme which is Asymptotic Preserving inside the domain and which provides a good approximation of the exact value at the boundary. To this purpose, we proceed as follows. Inside the space domain, we use the same numerical scheme as the one described in the previous section. At the boundary, we propose a specific treatment, without using neither artificial boundary condition nor any mesh refinement. 

Our approach is based on the following idea. We consider the spatial derivative of the density at the boundary as an additional unknown in the scheme. It is known indeed that the density becomes stiff (its derivative is of order $\frac{1}{\eps}$ in the boundary layer) and the analysis of the boundary layer requires usually a rescaling process. Here, our strategy enables to avoid such rescaling. 
In the boundary cells, we simply determine the slope of $\overline \rho$ by ensuring the balance between the incoming and outgoing half-fluxes of $f$. We recall indeed that, if $\chi$ is a solution of the half-space problem \fref{milne}, then one has $\pa_y\langle v\chi\rangle_V=0$ and then, from \fref{exact}, we deduce that $\langle v\chi\rangle_V=0$ and then, at $y=0$, we have $\int_{V_+(0)}v\chi d\mu=-\int_{V_-(0)}vf_b d\mu$. We shall use this idea in the time evolution context.

\bs
Let us now describe our numerical scheme. We assume that we know $\rho^n_i$ and $\overline\rho^n_i$ for $i=0,\ldots,N+1$, and $g^n_{i+1/2}$ for $i=0,\ldots N$. We also assume that we know the approximation of the function $g$ at the boundary, i.e. $g_0^{n}$ and $g_{N+1}^{n}$ for all $v\in V$. We will now show how to compute all these quantities at time $t_{n+1}$.

\bs
\ni
{\em First step: determination of the quantities at time iteration $n+1$, near the boundaries: $\rho_0^{n+1}$, $g^{n+1}_{1/2}$, $g^{n+1}_{0}$ and $\rho_{N+1}^{n+1}$, $g^{n+1}_{N+1/2}$, $g^{n+1}_{N+1}$}
\nopagebreak

\ms
\ni
First, we use a centered scheme for the calculation of $g^{n+1}_{1/2}$ and $g^{n+1}_{N+1/2}$:
\bea
\label{gbord1}
\frac{g_{1/2}^{n+1}-g_{1/2}^n}{\Delta t}&+&\frac{1}{\eps}(I-\Pi_{V_-})\left(v\frac{g_1^n-g_0^n}{\Delta x}\right)\\
&+&\frac{1}{\eps}(I-\Pi_{V_-})(v\calE)\left(\frac{\overline \rho^n_{1}-\overline \rho_{0}}{\Delta x}\right)=\frac{1}{\eps^2}(I-\Pi_{V_-})(Lg_{1/2}^{n+1}),
\nonumber\eea
\bea
\label{gbord2}
\frac{g_{N+1/2}^{n+1}-g_{N+1/2}^n}{\Delta t}&+&\frac{1}{\eps}(I-\Pi_{V_-})\left(v\frac{g_{N+1}^n-g_{N}^n}{\Delta x}\right)\\\nonumber
&+&\frac{1}{\eps}(I-\Pi_{V_-})(v\calE)\left(\frac{\overline \rho_{N+1}-\overline \rho^n_{N}}{\Delta x}\right)=\frac{1}{\eps^2}(I-\Pi_{V_-})(Lg_{N+1/2}^{n+1}),
\eea
where $g_1^n$ and $g_N^n$ are determined by interpolation
$$g_1^n=\frac{g_{1/2}^n+g_{3/2}^n}{2},\qquad g_{N}^n=\frac{g_{N-1/2}^n+g_{N+1/2}^n}{2}.$$
It remains to determine $g^{n+1}_{0}$ and $g^{n+1}_{N+1}$. For incoming velocities, these functions are prescribed thanks to \fref{inc}:
\bea
\label{datainc1}\un_{V_-(0)}\,g^{n+1}_0(v)&=&\un_{V_-(0)}\left(f_\ell(v)-\overline \rho_0\calE\right),\\
\label{datainc2}\un_{V_-(1)}\,g^{n+1}_{N+1}(v)&=&\un_{V_-(1)}\left(f_r(v)-\overline \rho_{N+1}\calE\right),
\eea
where $\un_A$ is the characteristic function of a set $A$. For outgoing velocities, we discretize the transport equation \fref{MM-g} in a suitable way. The specific point here is that we consider the derivative of $\overline \rho$ at the boundary as an unknown. We approximate $\eps \frac{\pa \overline \rho^n}{\pa x}(0)$ by the quantity $\lambda_\ell^n$ and $\eps \frac{\pa \overline \rho^n}{\pa x}(L)$ by the quantity $\lambda_r^n$, which are computed below. For $v<0$ only, we write
\bea
\hspace*{-1cm}\frac{g^{n+1}_0-g^n_0}{\Delta t}&+&\frac{\lambda_\ell^{n+1}}{\eps^2}(I-\Pi_{V_-})(v\calE)
\nonumber\\
&+&\frac{1}{\eps}(I-\Pi_{V_-})\left(v\frac{g^{n}_{1/2}-g^{n}_{-1/2}}{\Delta x}\right)=\frac{1}{\eps^2}(I-\Pi_{V_-})(Lg^{n+1}_0)
\label{bound1}
\eea
and, for $v>0$ only, we write
\bea
\hspace*{-1cm}
\frac{g^{n+1}_{N+1}-g^n_{N+1}}{\Delta t}&+&\frac{\lambda_r^{n+1}}{\eps^2}(I-\Pi_{V_-})(v\calE)\nonumber\\
&+&\frac{1}{\eps}(I-\Pi_{V_-})\left(v\frac{g^{n}_{N+3/2}-g^{n}_{N-1/2}}{\Delta x}\right)=\frac{1}{\eps^2}(I-\Pi_{V_-})(Lg^{n+1}_{N+1}).
\label{bound2}
\eea
In these equations, we have used the quantities $g_{-1/2}^{n}$ and $g_{N+3/2}^{n}$ which are again defined by extrapolation:
\be
\label{hhhj}
g_{-1/2}^{n}=2g_0^{n}-g_{1/2}^{n}\quad \mbox{and}\quad g_{N+3/2}^{n}=2g_{N+1}^{n}-g_{N+1/2}^{n}\,.
\ee
Let us give a sense to these two equations \fref{bound1} and \fref{bound2} which are both written on $V_+$. We only deal with the case of \fref{bound1}, the other equation can be treated similarly. This equation takes the form
$$\un_{V_+}\left(\frac{\eps^2}{\Delta t}I -\widetilde L\right)g^{n+1}_0=\un_{V_+}h,\quad \mbox{with}\quad \widetilde L=(I-\Pi_{V_-})L(I-\Pi_{V_-}).$$
We recall that $g^{n+1}_0$ is already known on $V_-$ by \fref{datainc1}. Hence, we have to solve the following linear system on the unknown $\un_{V_+} g^{n+1}_0$:
$$\un_{V_+}\left(\frac{\eps^2}{\Delta t}I -\widetilde L\right)\un_{V_+}g^{n+1}_0=\un_{V_+}h+\un_{V_+}\widetilde L\un_{V_-}g^{n+1}_0,$$
where the right-hand side is given. To see that this system is well posed, it suffices to recall that the operator $\widetilde L$ being self-adjoint and nonpositive, the bilinear form associated to $\frac{\eps^2}{\Delta t} I-\widetilde L$ is symmetric positive definite on the set of functions defined on $V$. Hence this bilinear form is also symmetric positive definite on the set of functions defined on $V_+$, and consequently the operator $\un_{V_+}\left(\frac{\eps^2}{\Delta t} -\widetilde L\right)\un_{V_+}$ is invertible on this set: Eq. \fref{bound1} admits a unique solution $\un_{V_+} g^{n+1}_0$.
From \fref{bound1}, we get
\bea
\un_{V_+}g^{n+1}_0&=&\left(\un_{V_+}\left(\frac{\eps^2}{\Delta t}I-\widetilde L\right)\un_{V_+}\right)^{-1}\left[\vphantom{\left(v\frac{g^{n}_{1/2}-g^{n}_{-1/2}}{\Delta x}\right)}\frac{\eps^2}{\Delta t}g^n_0-\lambda_\ell^{n+1}(I-\Pi_{V_-})(v\calE)\right.\nonumber\\&&\hspace*{1.5cm}\left.-\eps(I-\Pi_{V_-})\left(v\frac{g^{n}_{1/2}-g^{n}_{-1/2}}{\Delta x}\right)+\widetilde L\un_{V_-}g^{n+1}_0\right].\label{inverse}
\eea
In order to determine $\lambda^{n+1}_\ell$ and $\lambda^{n+1}_r$, we discretize \fref{MM-rho} at $x=0$ and at $x=1$: 
$$
\frac{\rho_0^{n+1}-\rho_0^{n}}{\Delta t}+\frac{1}{\eps}\left\langle v\frac{g^{n+1}_{1/2}-g^{n+1}_{-1/2}}{\Delta x}\right\rangle_V=0,
$$
$$
\frac{\rho_{N+1}^{n+1}-\rho_{N+1}^{n}}{\Delta t}+\frac{1}{\eps}\left\langle v\frac{g^{n+1}_{N+3/2}-g^{n+1}_{N+1/2}}{\Delta x}\right\rangle_V=0.
$$
Using \fref{rhobarrho}, \fref{incbar} and \fref{hhhj}, we deduce the following two equations:
\be
\label{bound3}
\left\langle \frac{g_0^{n+1}-g_0^{n}}{\Delta t}\right\rangle_V+\frac{2}{\eps}\left\langle v\frac{g^{n+1}_{1/2}-g^{n+1}_{0}}{\Delta x}\right\rangle_V=0
\ee
\be
\label{bound4}
\left\langle \frac{g_{N+1}^{n+1}-g_{N+1}^{n}}{\Delta t}\right\rangle_V+\frac{2}{\eps}\left\langle v\frac{g^{n+1}_{N+1}-g^{n+1}_{N+1/2}}{\Delta x}\right\rangle_V=0.
\ee
Inserting \fref{inverse}, together with \fref{datainc1}, into \fref{bound3} enables to compute explicitely $\lambda_\ell^{n+1}$, from which we get  $g^{n+1}_0$ using \fref{inverse}. Finally, we compute $\rho_0^{n+1}$ using \fref{rhobarrho} and \fref{incbar}. We proceed similarly with \fref{bound2} and \fref{bound4} to compute $\lambda_r^{n+1}$, $g^{n+1}_{N+1}$ and $\rho^{n+1}_{N+1}$.

\bs
\ni
{\em Second step: determination of the quantities at time iteration $n+1$, inside the domain.}
\nopagebreak

\ms
\ni
In this second step, we focus on the interior of the domain and compute the values $\rho^{n+1}_i$, $\overline\rho^{n+1}_i$ for $i=1, \ldots, N$ and $g^n_{i+1/2}$ for $i=1, \ldots, N-1$.

For $i=1,\ldots,N-1$, we write
\bea
\hspace*{-10mm}\frac{g_{i+1/2}^{n+1}-g_{i+1/2}^n}{\Delta t}&+&\frac{1}{\eps}(I-\Pi_{V_-})\left(v^+\frac{g_{i+1/2}^n-g_{i-1/2}^n}{\Delta x}+v^-\frac{g_{i+3/2}^n-g_{i+1/2}^n}{\Delta x}\right)\nonumber\\
&+&\frac{1}{\eps}(I-\Pi_{V_-})(v\calE)\left(\frac{\overline \rho^n_{i+1}-\overline \rho^n_{i}}{\Delta x}\right)=\frac{1}{\eps^2}(I-\Pi_{V_-})(Lg_{i+1/2}^{n+1}).
\label{ginter}
\eea
Then, for $i=1,\ldots,N$, we write 
\be
\label{rhointer}
\frac{\rho^{n+1}_i-\rho^n_i}{\Delta t}+\frac{1}{\eps}\cint{v\frac{g_{i+1/2}^{n+1}-g_{i-1/2}^{n+1}}{\Delta x}}_V=0.
\ee
Note that the values $g^{n+1}_{1/2}$ and $g^{n+1}_{N+1/2}$ have been computed in the first step.

\begin{figure}[h]
\begin{center}
\includegraphics[width=12cm]{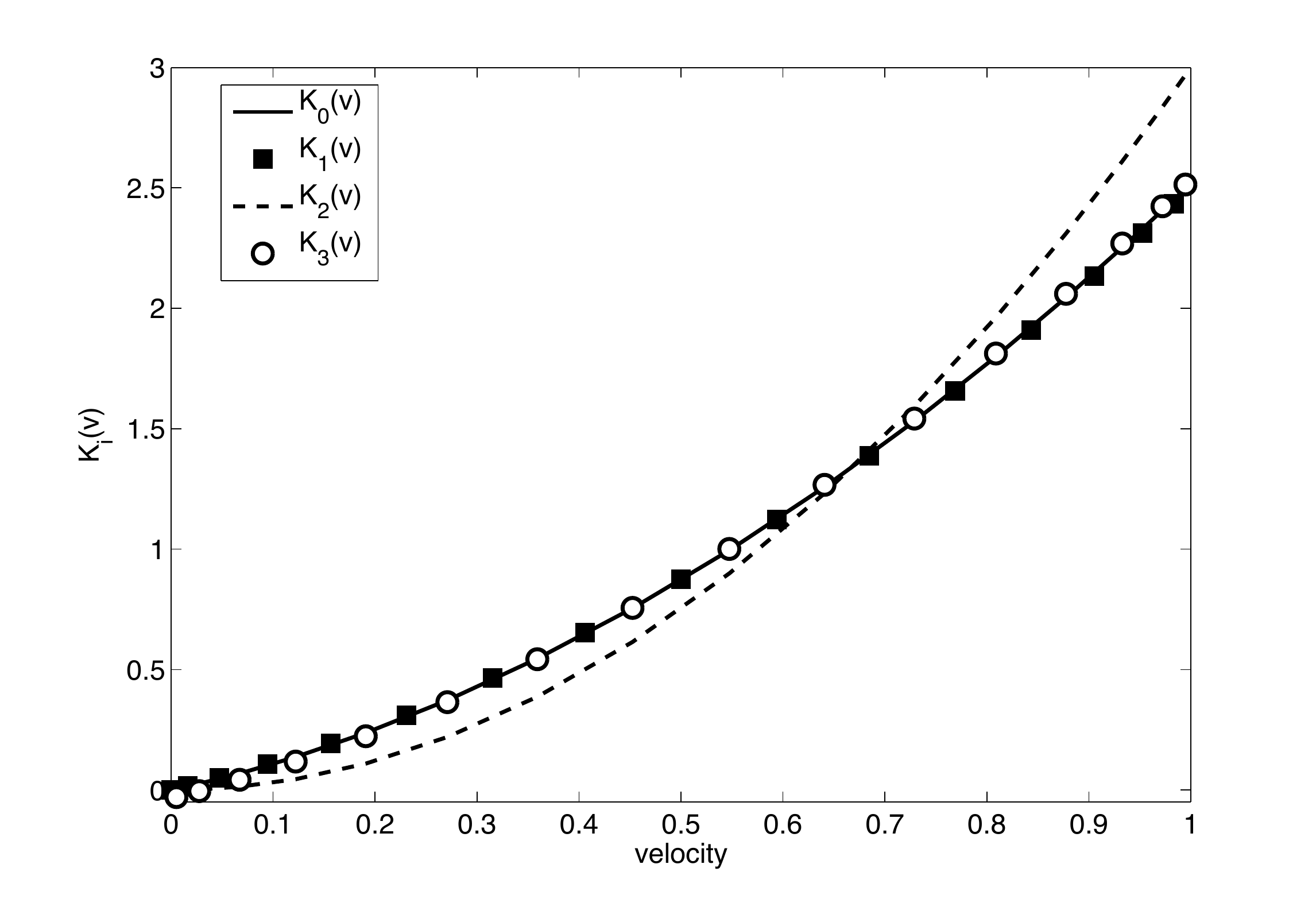}
\caption{Exact kernel $K_0(v)$ given by \fref{chandra} and its approximations $K_1$, $K_2$, $K_ 3$ defined respectively by \fref{polynome1}, \fref{polynome2}, \fref{polynome3}}\label{fig1}
\end{center}
\end{figure}

\begin{proposition}[Formal]
\label{propschema2}
The numerical scheme \fref{gbord1}, \fref{gbord2}, \fref{bound1}, \fref{bound2}, \fref{bound3}, \fref{bound4}, \fref{ginter}, \fref{rhointer} with initial and boundary conditions \fref{initialdisc}, \fref{incbar}, \fref{inc} determines all the values $\rho^n_i$ and  $\overline\rho^n_i$ for $i=0,\ldots,N+1$, and $g^n_{i+1/2}$ for $i=0,\ldots N$, at any time interation $n$. This scheme is stable under the CFL condition
$$\Delta t\leq C\left(\Delta x^2+\eps\Delta x\right),$$
where $C$ is a constant independent of $\eps$, $\Delta t$ and $\Delta x$. For all fixed $\eps>0$, this scheme is consistent with the initial and boundary value problem \fref{eq-transp}, \fref{initialdata}, \fref{BC-f}. Moreover, for all fixed $\Delta t$ and $\Delta x$, the scheme degenerates as $\eps \to 0$ into a discretization of the diffusion equation \fref{eq-diff}. More precisely, for $i=1,\ldots, N$, we have  $\rho_i^n\to  \widehat \rho_i^n$ as $\eps \to 0$, and $\widehat \rho^n_i$ is defined by the scheme
$$
\frac{\widehat \rho^{n+1}_i-\widehat \rho^n_i}{\Delta t}-\kappa\frac{\widehat \rho^n_{i+1}-2\widehat \rho^n_i+\widehat \rho^n_{i-1}}{\Delta x^2}=0,\quad \mbox{for}\quad i=1,\ldots,N,
$$
where $\kappa>0$ is given in \fref{eq-diff} and where the limiting boundary values take the form
\be
\label{eqlimite0}
\widehat \rho_0^n=\int_{V_-(0)}K_L(v)f_\ell(v)dv, \qquad \widehat \rho_{N+1}^n=\int_{V_-(1)}K_L(-v)f_r(v)dv,
\ee
where $K_L$ is an explicit function which depends only on the operator $L$. The expression of this kernel can be deduced from \fref{lambdafinal} and \fref{eqfin} below.
\end{proposition}
In the particular case where $V=[-1,1]$, $\calE=1$, $d\mu=\frac{1}{2}dv$ and $L=\Pi-I$, we obtain $\kappa=\frac{1}{3}$ and 
we can compute the following expressions for the kernels:
$$
\widehat \rho_0^n=\int_0^1K_3(v)f_\ell(v)dv,\quad \widehat \rho_{N+1}^n=\int_{-1}^0K_3(-v)f_r(v)dv.
$$
with
\be
\label{polynome3}
K_3(v)=\frac{3}{2}v^2+\frac{15}{14}v-\frac{1}{28}.
\ee
For instance, if $f_\ell(v)=v$, we get the value $\widetilde \rho_0=0.7143$ as boundary value for the limiting diffusion equation, where we recall that the right value is given by $\rho_b(0)=\frac{\sqrt{3}}{2}\int_0^1 v^2H(v)dv\simeq 0.7104$. Therefore, this scheme provides a much better approximation of the limiting boundary value than the one obtained from the first scheme presented in Section \ref{sectattempt}. On Figure \ref{fig1}, we plotted four curves: the exact kernel $K_0(v)$ defined by \fref{chandra}, its usual approximation $K_1(v)$ defined by \fref{polynome1}, its approximation $K_2(v)$ defined by \fref{polynome2} and given by the first scheme presented in Section \ref{sectattempt}, and its approximation $K_3(v)$ defined by \fref{polynome3}. We observe that our kernel $K_3(v)$, as well as $K_1(v)$, fit very well with the exact Chandrasekhar kernel $K_0(v)$, whereas the kernel $K_2(v)$ does not provide a good approximation.

\bs
\ni
{\em Proof of Proposition \pref{propschema2}.}
The proof is similar to that of Proposition \ref{propschema1}, except that we have to deal with the boundary terms. We focus on the left boundary $x=0$, the right boundary can be treated similarly. Let $\eps>0$ be fixed. Our numerical scheme is consistent with the continuous model as soon as the following property holds: 
\be
\label{identlambda}
\lambda_\ell^{n}=\eps \frac{\overline \rho_1^n-\overline \rho_0}{\Delta x}+\calO(\Delta t,\Delta x).
\ee
To show this property, we consider a modified scheme consisting in the same equations \fref{gbord1}, \fref{gbord2}, \fref{ginter}, \fref{rhointer}, and where \fref{bound1} is replaced by
\bea
\hspace*{-1cm}\frac{g^{n+1}_0-g^n_0}{\Delta t}&+&\frac{1}{\eps}(I-\Pi_{V_-})(v\calE)\left(\frac{\overline \rho_1^{n+1}-\overline \rho_0}{\Delta x}\right)
\nonumber\\
&+&\frac{1}{\eps}(I-\Pi_{V_-})\left(v\frac{g^{n}_{1/2}-g^{n}_{-1/2}}{\Delta x}\right)=\frac{1}{\eps^2}(I-\Pi_{V_-})(Lg^{n+1}_0)
\label{bound1modif}
\eea
for $v<0$. This equation allows to compute $g^{n+1}_0(v)$ for $v<0$. Since it is clearly a consistent discretization of \fref{MM-g}, the whole modified scheme will be consistent with the continuous model and therefore satisfies
$$\left\langle \frac{g_0^{n+1}-g_0^{n}}{\Delta t}\right\rangle_V+\frac{2}{\eps}\left\langle v\frac{g^{n+1}_{1/2}-g^{n+1}_{0}}{\Delta x}\right\rangle_V=\calO(\Delta t,\Delta x).$$
Now, let us recall that, in our numerical scheme, the quantity $\lambda_\ell^{n}$ was computed thanks to \fref{bound1} and \fref{bound3}. This enables to conclude that $\lambda_\ell^{n+1}$ is a discretization of $\eps \frac{\pa \overline \rho^{n+1}}{\pa x}$ and \fref{identlambda} holds.

Let us now analyze the limit $\eps\to 0$ of our scheme near the left boundary. From \fref{gbord1}, we get
\be
\label{g12}
g_{1/2}^{n+1}=\eps(I-\Pi_{V_-})L^{-1}\left[v\calE\left(\frac{\overline \rho^n_{1}-\overline \rho_{0}}{\Delta x}\right)-\frac{1}{\Delta x}(I-\Pi)(vg_0^n)\right]+\calO(\eps^2).
\ee
Now, we need to pass to the limit in \fref{inverse}. For this, we recall that the bilinear form associated to $-\widetilde L=-(I-\Pi_{V_-})L(I-\Pi_{V_-})$ is symmetric and positive definite on the set $\{f:\,\langle f\rangle_{V_-}=0\}$. Hence, since the set of functions which vanish on $V_-$ is included in this set, the operator $\un_{V_+}\widetilde L\un_{V_+}$ is invertible from the set of functions defined on $V_+$ to itself. Thus, from \fref{datainc1} and \fref{inverse}, we obtain
\be
\label{jkl}
\hspace*{-5mm}\begin{array}{lll}
\un_{V_-}g_0^{n}(v)&=&\un_{V_-}\left(f_\ell(v)-\overline \rho_0\calE\right),\\[3mm]
\un_{V_+}g_0^{n}(v)&=&\left(\un_{V_+}\widetilde L\un_{V_+}\right)^{-1}\left[\lambda_\ell^{n}(I-\Pi_{V_-})(v\calE)-\widetilde L\un_{V_-}g^{n}_0\right]+\calO(\eps).
\end{array}
\ee
It remains to determine $\lambda_\ell^{n}$. From \fref{bound3}, we get the equality of the half-fluxes at the boundary:
$$\left\langle v g_0^{n}\right\rangle_V=\calO(\eps).$$
Inserting \fref{jkl} into this equation yields $\lambda^{n}_\ell=\lambda_\ell+\calO(\eps)$, where
\be
\label{lambdafinal}\lambda_\ell=\frac{\left\langle -v^+(f_\ell -\overline \rho_0\calE)+v^-(\un_{V_+}\widetilde L\un_{V_+})^{-1}\left[\widetilde L(\un_{V_-}(f_\ell -\overline \rho_0\calE))\right]\right\rangle_V}{\left\langle v^-(\un_{V_+}\widetilde L\un_{V_+})^{-1}\left[(I-\Pi_{V_-})(v\calE)\right]\right\rangle_V}.
\ee
Next, equation \fref{rhointer} on $\rho_1^{n+1}$ becomes
$$
\frac{\rho^{n+1}_1-\rho^n_1}{\Delta t}-\kappa\frac{\rho^n_{2}-2\rho^n_1+\overline\rho^n_0}{\Delta x^2}-\frac{\kappa}{\Delta x^2}\,\frac{\left\langle vL^{-1}(I-\Pi)(vg_0^n)\right\rangle_V}{\left\langle vL^{-1}(v\calE)\right\rangle_V}=\calO(\eps)
$$
where we used the expression \fref{g32} for $g^{n+1}_{3/2}$ and \fref{g12} for $g^{n+1}_{1/2}$.
It remains to replace $g_0^n$ in this formula by its expression \fref{jkl}, in which $\lambda^n_\ell$ is given by \fref{lambdafinal}. We finally get
$$
\frac{\rho^{n+1}_1-\rho^n_1}{\Delta t}-\kappa\frac{\rho^n_{2}-2\rho^n_1+\widehat\rho^n_0}{\Delta x^2}=\calO(\eps),
$$
where $\widehat\rho^n_0$ is given by
\bea
\label{eqfin}
\widehat \rho_0^n&=&\overline\rho_0-\frac{1}{\kappa}\left\langle vL^{-1}(I-\Pi)\left[v^+(f_\ell-\overline\rho_0\calE)\right]\right\rangle_V\\&&\hspace*{-5mm}-\frac{1}{\kappa}\left\langle vL^{-1}(I-\Pi)\left[ v^-\left(\un_{V_+}\widetilde L\un_{V_+}\right)^{-1}\left(\lambda_\ell(I-\Pi_{V_-})(v\calE)-\widetilde L\un_{V_-}(f_\ell-\overline\rho_0\calE)\right)\right]\right\rangle_V.\nonumber
\eea
Finally, we remark that inserting the expression \fref{lambdafinal} of $\lambda_\ell$ into \fref{eqfin} yields the existence of the kernel $K_L$ such that \fref{eqlimite0} holds. The proof of Proposition \ref{propschema2} is complete. 
\qed

\bs
Let us compute the approximate boundary value given by our scheme in a special case of interest, which is tested numerically below. We set $V=[-1,1]$, $d\mu=\frac12 dv$ and take the collision operator under the form
$$Lf(v)=\frac{1}{2}\int_{-1}^1\sigma(v,w)(f(w)-f(v))dw\quad \mbox{with }\sigma(v,w)=p(v)p(w),$$
$p$ being a given even function. Then straightforward but lengthy computations from \fref{lambdafinal} and \fref{eqfin} yield
$$
\widehat \rho_0^n=\int_0^1K_\sigma^{app}(v)f_\ell(v)dv,
$$
where the kernel is given by
\be
\label{polynomesigma2}
K^{app}_\sigma(v)=a\frac{v^2}{p(v)}+bv+cp(v)+d.
\ee
The coefficients $a$, $b$, $c$  and $d$ are given by
$$a=\frac{1}{2\kappa\alpha},\quad b=\frac{\gamma+2\alpha \delta}{\kappa\alpha(1+4\alpha\gamma)},\quad c=\frac{2\gamma^2-\delta}{\kappa\alpha(1+4\alpha\gamma)},\quad d=1-\frac{\gamma}{\kappa \alpha},$$
where we have set
$$\alpha=\int_0^1p(v)dv,\quad \beta=\int_0^1\frac{v}{p(v)}dv,\quad \gamma=\int_0^1\frac{v^2}{p(v)}dv,\quad \delta=\int_0^1\frac{v^3}{(p(v))^2}dv.$$
Notice that in the simplified case where $L=\Pi-I$, i.e. $p(v)\equiv 1$, then one has $a=\frac{3}{2}$, $b=\frac{15}{14}$, $c=-\frac{1}{28}$ and $d=0$, which covers the expression \fref{polynome3}.

\section{Numerical tests}
\label{sectionnum}

We present here some numerical experiments to validate the approach. In our tests, the space domain is the interval $\Omega=[0,1]$, the velocity domain is $V=[-1,1]$ and $d\mu=\frac12 dv$. All the integrals in velocity which are involved in the model are computed using a Gauss quadrature method. The collision operator will be of the following form:
$$Lf(v)=\frac{1}{2}\int_{-1}^1\sigma(v,w)(f(w)-f(v))dw,$$
where $\sigma$ is a given symmetric and nonnegative kernel. We have $\mathcal E(v)\equiv 1$. The initial data will be zero: $f_{init}(x,v)=0$. The function $\omega(x,v)$ used in our simulations is $\omega(x,v)=(2x-1)v-x(1-x)$. In the tests using our numerical scheme, the time step is linked to the space grid step as follows: $\Delta t=\frac{1}{2}(\frac{\Delta x^2}{2}+\eps \Delta x)$.

We compare our numerical results, which will be referred to as LMe, with those given by the following schemes:
\begin{enumerate}[i.]
\item The time explicit upwind scheme for the original kinetic equation \fref{eq-transp}, where $\Delta t$ and $\Delta x$ are linked to $\eps$ by standard CFL stability condition. When $\eps\geq 10^{-2}$, this explicit scheme will be highly resolved to serve as a reference for our comparisons.
\item The direct discretization of the diffusion equation \fref{eq-diff} by a standard explicit scheme, in the cases where the right boundary condition is known (see e.g. the Chandrasekhar value \fref{chandra1}). Our results will be compared with those given by this scheme when $\eps$ is very small.
\item The micro-macro scheme of Lemou and Mieussens given in \cite{LM-AP1}, which will be referred to as LMi.
\item The scheme by Klar \cite{Klar1,Klar3}, which will be referred to as K.
\item The scheme by Jin, Pareschi and Toscani \cite{JPT}, which will be referred to as JPT.
\end{enumerate}
In all the figures, we plot the density $\rho$ as a function of space, at different times.

\subsection{Boundary conditions at thermal equilibrium}

In this subsection, we present some numerical tests in a situation where the incoming boundary condition is proportional to the equilibrium $\calE(v)\equiv 1$. We recall that, in this case, there is no boundary layer in the diffusive limit.

\bs
\ni
{\em Example 1.} Kinetic regime with isotropic boundary conditions:
$$
\sigma(v,w)=1,\quad f_\ell(v)=1,\quad f_r(v)=0,\quad \eps=1.
$$
On Figure \ref{fig2}, we plot the results obtained with our scheme LMe and the reference obtained by the explicit scheme at times $t=0.2$, $t=0.5$, $t=1$, $t=2$ and $t=4$. Our scheme LMe is simulated with 50 and 200 grid points in space, whereas the explicit scheme is simulated with 2000 grid points. We clearly see that our results are in good agreement with the reference simulation in this kinetic regime.

\bs
\ni
{\em Example 2.} Diffusive regime with isotropic boundary conditions:
$$
\sigma(v,w)=1,\quad f_\ell(v)=1,\quad f_r(v)=0,\quad \eps=10^{-4}.
$$
On Figure \ref{fig3}, we plot the results obtained with our scheme LMe and the reference obtained by the diffusion equation at times $t=0.01$, $t=0.1$ and $t=0.4$. Our scheme LMe is simulated with 50 and 200 grid points in space, whereas the diffusion limit is simulated with 200 grid points. We also see here that our results are in total agreement with the diffusion limit. Of course, the AP property of our scheme allows us to keep $\Delta t$ and $\Delta x$ independent of $\eps$.

\subsection{Non equilibrium boundary conditions and boundary layers}
\label{sigmaconstant}
In this subsection, we present some numerical tests where the boundary data are not at equilibrium, which induce a boundary layer in the diffusive regimes. Here, we still consider the situation where the cross section $\sigma(v,w)$ is constant, we differ to the next section the more general case. 

\bs
\ni
{\em Example 3.} Kinetic regime with non isotropic boundary conditions:
$$
\sigma(v,w)=1,\quad f_\ell(v)=v,\quad f_r(v)=0,\quad \eps=1.
$$
On Figure \ref{fig4}, we plot the results obtained with our scheme LMe and the reference solution obtained by the explicit scheme at times $t=0.2$, $t=0.4$ and $t=0.8$. Our scheme LMe is simulated with 50 and 200 grid points in space, whereas the explicit scheme is simulated with 2000 grid points. Again, we observe a good agreement between the two methods, in the kinetic regime.

\bs
\ni
{\em Example 4.} Intermediate regime with non isotropic boundary conditions:
$$
\sigma(v,w)=1,\quad f_\ell(v)=v,\quad f_r(v)=0,\quad \eps=10^{-2}.
$$
On Figures \ref{fig5} and \ref{fig6}, we plot the results obtained with our scheme LMe, the schemes K, JPT, LMi with 10 gridpoints or 50 gridpoints in space, and compare them to the reference solution obtained with the explicit scheme using 10000 gridpoints in space, at time $t=0.4$. It is clear that our scheme is in a good agreement with the reference, even with coarse grids. There is no need in our scheme to discretize the boundary layer. Moreover, we recall that no artificial boundary condition is imposed, since our scheme is constructed in such a way that only the natural inflow boundary condition is used. We also see on Figures \ref{fig5} and \ref{fig6} that the numerical boundary value (at $x=0$) accurately approximates the right value and we have the following remarkable property: the curve obtained by our scheme is close to the reference curve whatever the mesh size is. This property is fulfilled in all our numerical tests.

\bs
\ni
{\em Example 5.} Diffusive regime with non isotropic boundary conditions:
$$
\sigma(v,w)=1,\quad f_\ell(v)=v,\quad f_r(v)=0,\quad \eps=10^{-4}.
$$
On Figures \ref{fig7} and \ref{fig8}, we plot the results obtained with our scheme LMe, the schemes K, JPT, LMi with 50 gridpoints or 200 gridpoints in space, and compare them to the reference solution obtained using a scheme for the diffusion limiting equation (with the exact Chandrasekhar value at the left boundary), at time $t=0.4$. Again, we observe that the results of our scheme LMe fit well with the reference curve. For the scheme K, we observe that the boundary value is very close to the right one, but the diffusive behavior is not correctly reproduced. For the schemes JPT and LMi, we observe that the uncorrectness of the diffusive regime is due to the fact that boundary condition is not the right one.

On Figure \ref{fig9}, we plot the density obtained by our scheme LMe using different numbers of space gridpoints: 10, 25, 50, 100 and 200, compared to the diffusion. We clearly observe that our scheme has a good behavior even with coarse grids.

\subsection{Collision operators with non constant cross sections}

In this subsection, we provide some numerical simulations in cases where the cross section $\sigma$ is non constant. Various expressions of $\sigma(v,w)$ will be tested. In the first two examples, the cross section takes a special form allowing to derive a formula similar to \fref{pointfixe} in order to compute numerically the associated generalized Chandrasekhar function. In these cases, we accurately know the exact value at the boundary in the diffusive regime, and this enables to compare our numerical results with this value. In the other two cases, we shall only compare our results with those obtained by the explicit scheme, in the intermediate regime $\eps=10^{-2}$. In all cases, we only consider a non isotropic boundary condition at $x=0$: $f_\ell(v)=v$, $f_r(v)=0$, which induces a boundary layer in the diffusive regime.

\bs
Before presenting the numerical tests, let us give the generalized Chandrasekhar like formula corresponding to the specific cases where the cross section has the form $\sigma(v,w)=p(v)p(w)$, $p$ being an even function. Following the method presented for instance in \cite{dautray-lions,golse1}, one can derive the following formula. Let us denote the limiting boundary value at $x=0$ for the diffusion problem by
$$\rho_\ell(0)=\int_0^1K_\sigma(v)f_\ell(v)dv.$$
Then the kernel $K_\sigma$ can be determined thanks to the generalized Chandrasekhar function $H_\sigma$ as follows:
\be
\label{Ksigma}
K_\sigma(v)=\frac{1}{2}vH_\sigma(v)\left(\int_0^1p(w)dw\right)^{-1/2}\left(\int_0^1\frac{w^2}{p(w)}dw\right)^{-1/2},
\ee
where $H_\sigma$ satisfies the nonlinear equation
\be
\label{chandra2}
H_\sigma(v)=1+\frac{1}{2\int_0^1p(w)dw}vH_\sigma(v)\int_0^1\frac{(p(w))^2}{vp(w)+wp(v)}H_\sigma(w)dw.
\ee
We shall compare numerically this kernel $K_\sigma(v)$ with its approximation $K^{app}_\sigma(v)$ given by our scheme, which has been computed above and is given by \fref{polynomesigma2}.

\bs
\ni
{\em Example 6:} $\sigma(v,w)=|v|^{3/2}|w|^{3/2}$.

\begin{figure}[h!]
\begin{center}
\includegraphics[width=12cm]{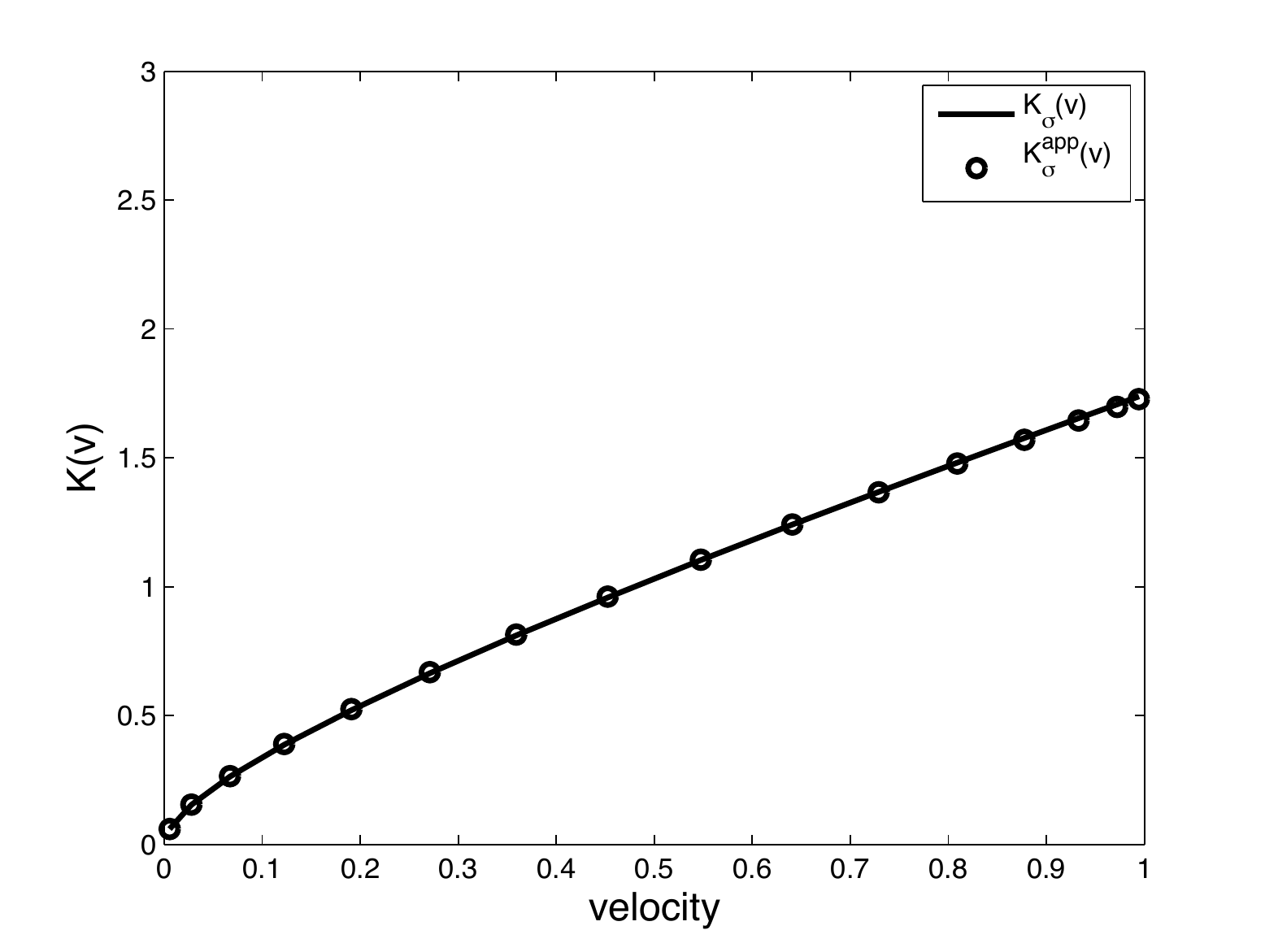}
\caption{Example 6. Exact kernel $K_\sigma(v)$ given by \fref{chandra2} and its approximation $K_\sigma^{app}$ defined by \fref{polynomesigma2}.}\label{figchandra2}
\end{center}
\end{figure}

\ms
\ni
We first plot on Figure \ref{figchandra2} the approximated kernel $K_\sigma^{app}(v)$ given by \fref{polynomesigma2}, compared to the exact kernel $K_\sigma(v)$ computed thanks to \fref{Ksigma}, \fref{chandra2}. We observe that the kernel generated by our scheme is a very accurate approximation of the exact kernel.

We plot the density in the same three different regimes as above: $\eps=1$ on Figure \ref{fig10}, $\eps=10^{-2}$ on Figure \ref{fig11} and $\eps=10^{-4}$ on Figure \ref{fig12}. In all cases, the observations we made in the Subsection \ref{sigmaconstant} for a constant collision kernel are still valid here and our curves are in good agreement with the reference.

\bs
\ni
{\em Example 7:} $\sigma(v,w)=(1-|v|^2)^{-3/4}(1-|w|^2)^{-3/4}$.

\begin{figure}[h!]
\begin{center}
\includegraphics[width=12cm]{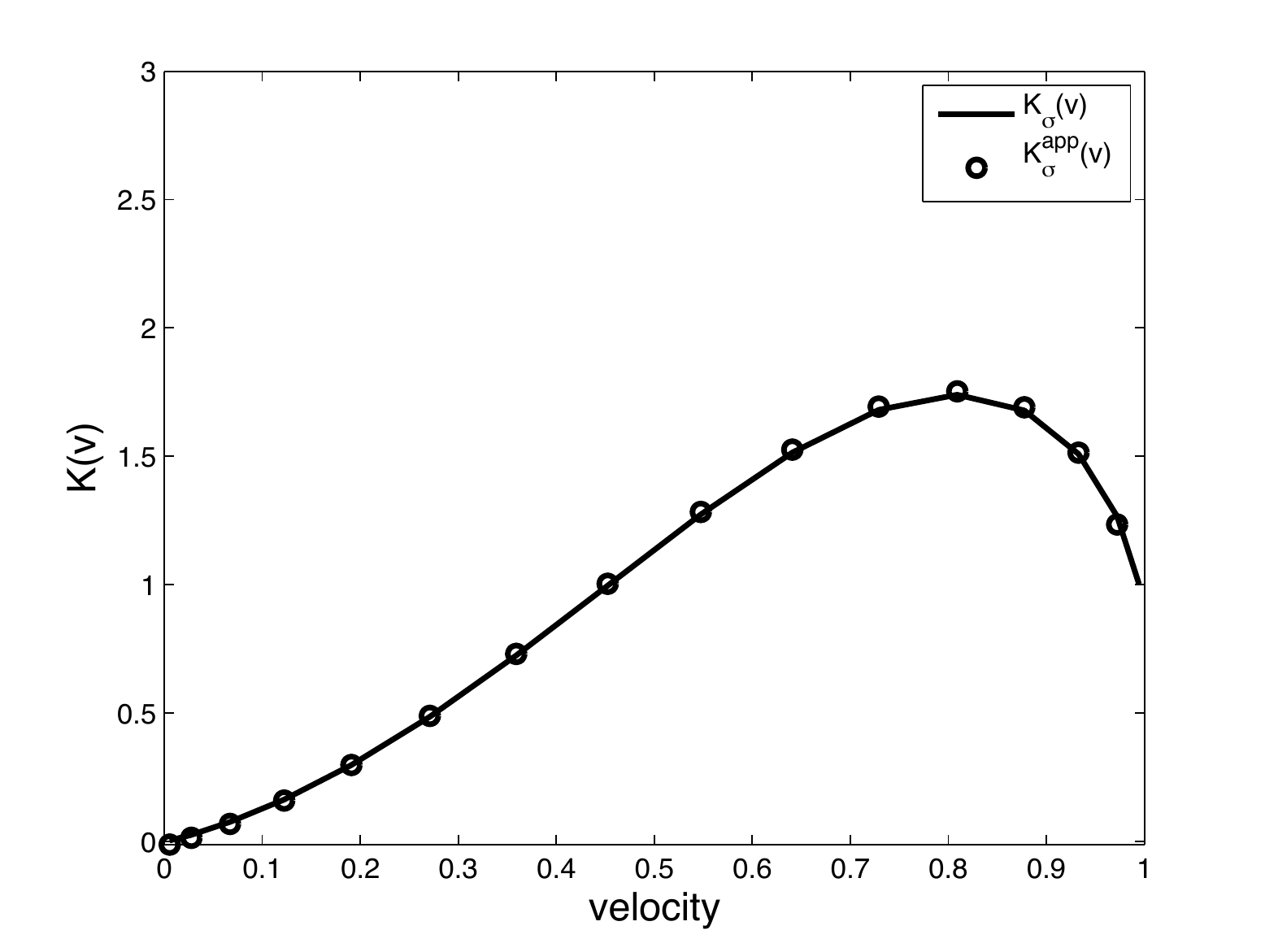}
\caption{Example 7. Exact kernel $K_\sigma(v)$ given by \fref{chandra2} and its approximation $K_\sigma^{app}$ defined by \fref{polynomesigma2}}\label{figchandra3}
\end{center}
\end{figure}

\begin{figure}[H!]
\begin{center}
\includegraphics[width=11cm]{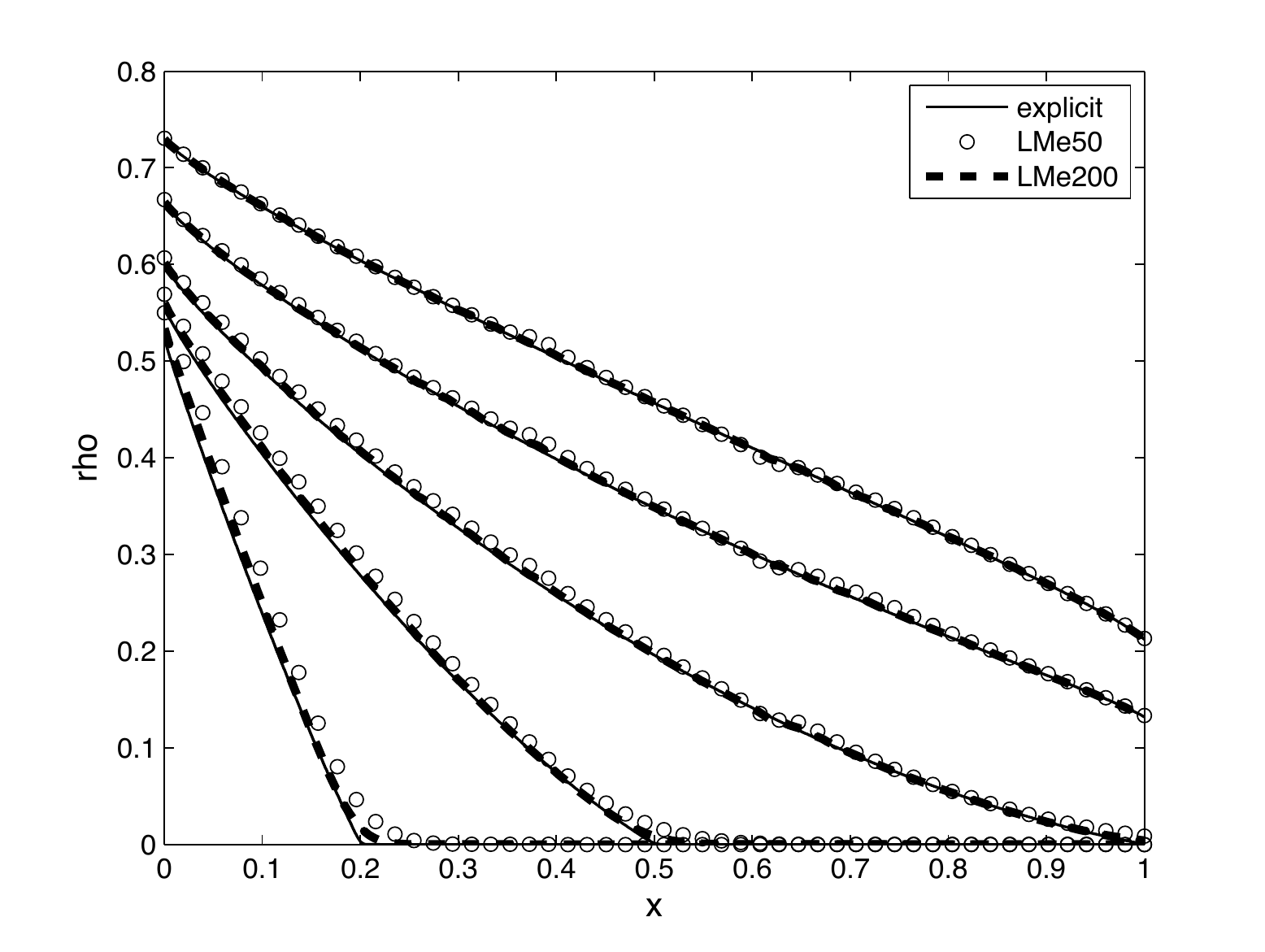}
\caption{Example 1: comparison between the explicit scheme and our scheme LMe (50 and 200 grid points), $\eps=1$. Results at times $t=0.2$, $t=0.5$, $t=1$, $t=2$ and $t=4$.}\label{fig2}
\end{center}
\end{figure}

\begin{figure}[H!]
\begin{center}
\includegraphics[width=11cm]{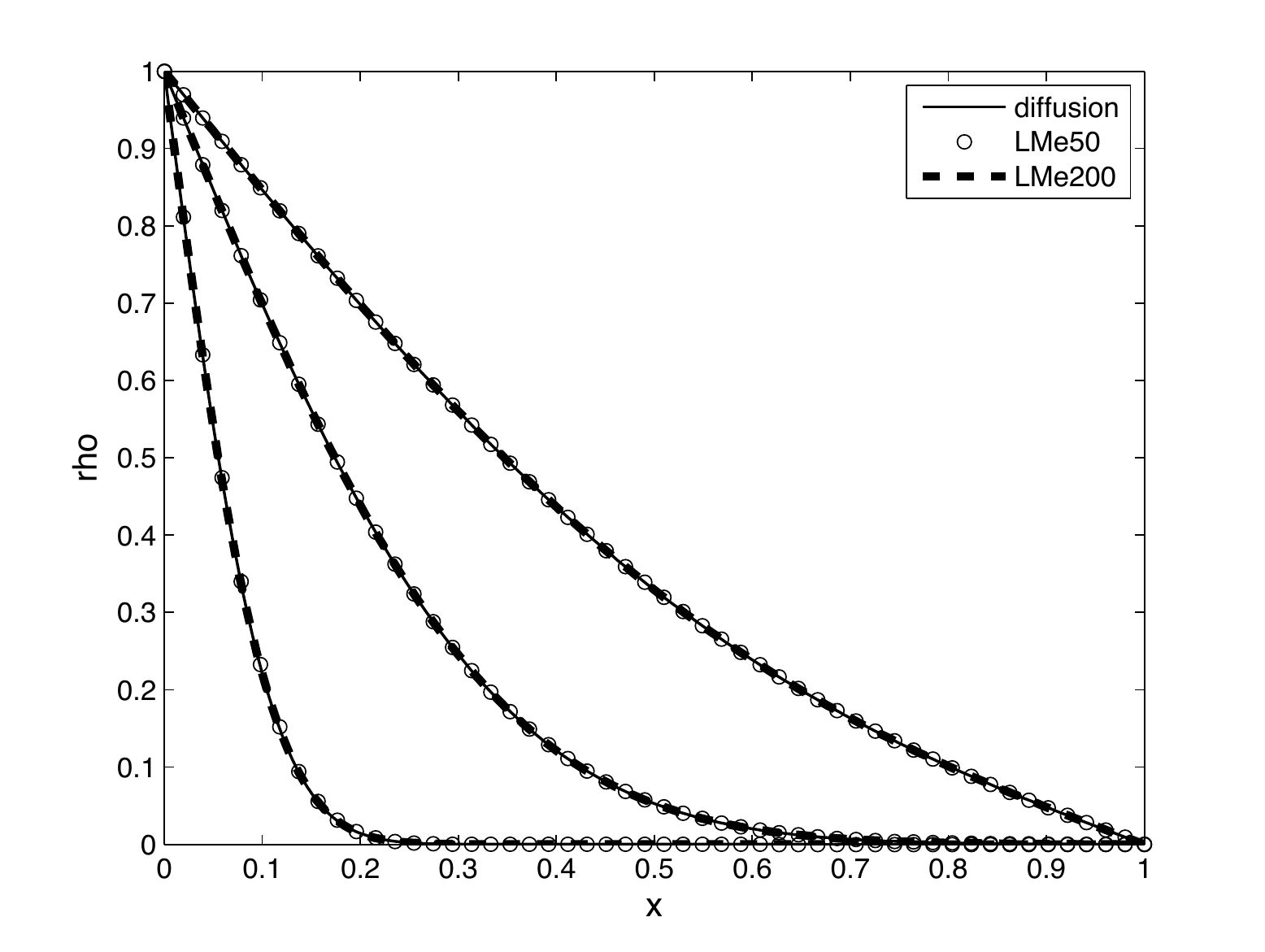}
\caption{Example 2: comparison between the diffusion limit and our scheme LMe (50 and 200 space gridpoints), $\eps=10^{-4}$. Results at times $t=0.01$, $t=0.1$ and $t=0.4$.}\label{fig3}
\end{center}
\end{figure}

\begin{figure}[H!]
\begin{center}
\includegraphics[width=11cm]{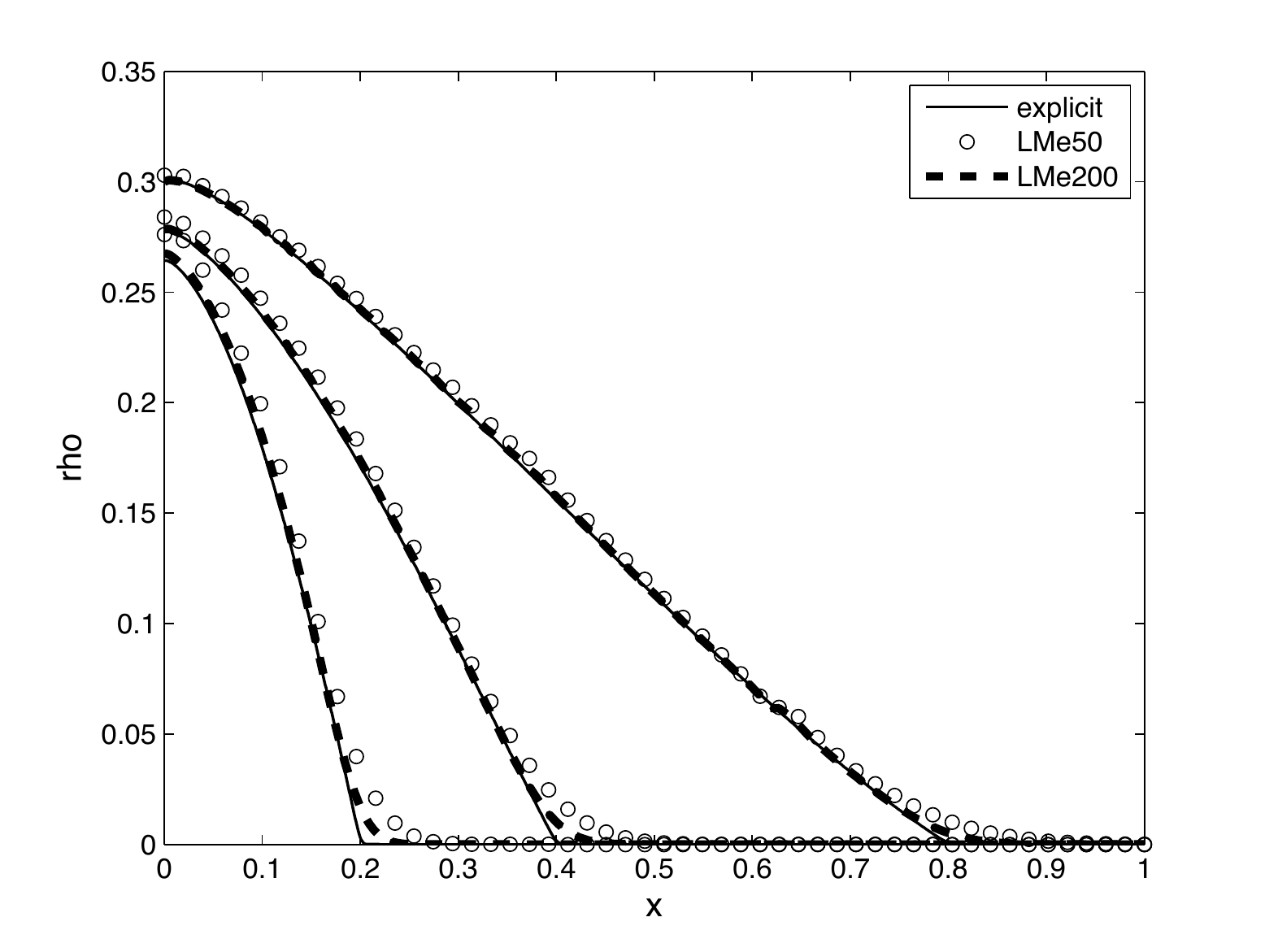}
\caption{Example 3: comparison between the explicit scheme and our scheme LMe (50 and 200 space gridpoints), $\eps=1$. Results at times $t=0.2$, $t=0.4$ and $t=0.8$.}\label{fig4}
\end{center}
\end{figure}

\begin{figure}[H!]
\begin{center}
\includegraphics[width=11cm]{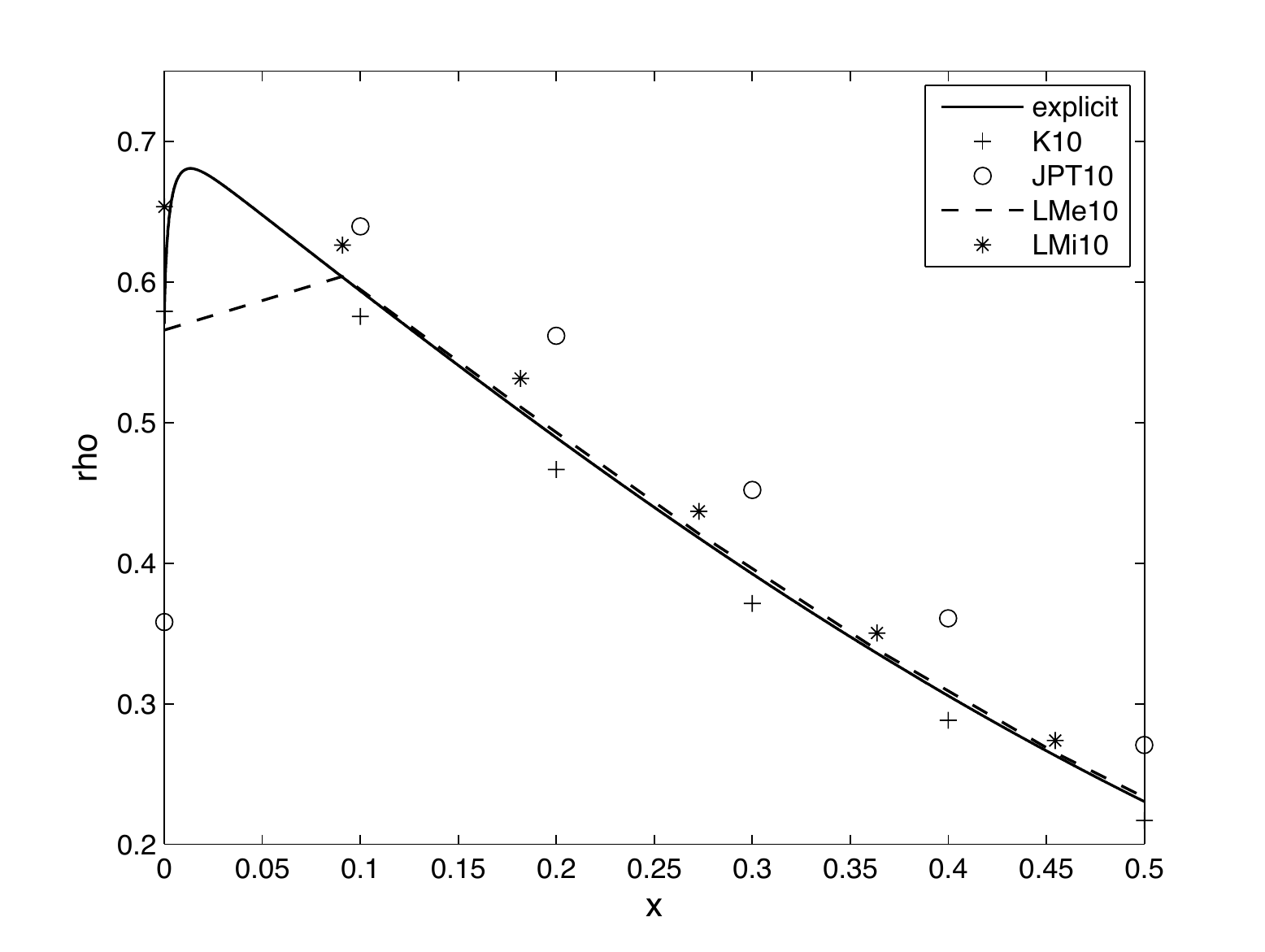}
\caption{Example 4: comparison between the explicit scheme (10 000 space gridpoints) and the schemes K, JPT, LMe and LMi (10 space gridpoints each), $\eps=10^{-2}$. A zoom has been made on the axes: $x\in [0,0.5]$.}\label{fig5}
\end{center}
\end{figure}

\begin{figure}[H!]
\begin{center}
\includegraphics[width=11cm]{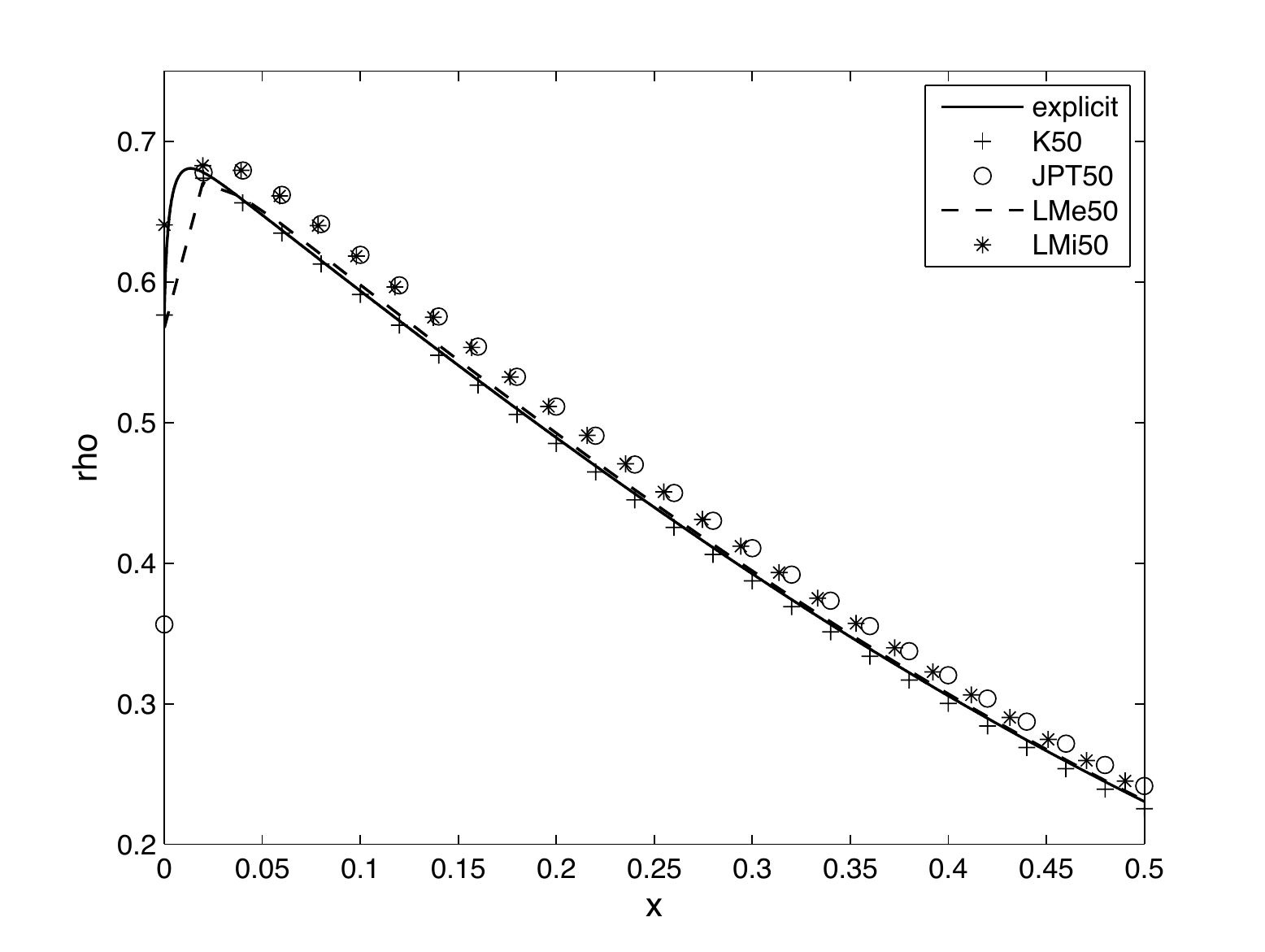}
\caption{Example 4: comparison between the explicit scheme (10 000 space gridpoints) and the schemes K, JPT, LMe and LMi (50 space gridpoints each), $\eps=10^{-2}$. A zoom has been made on the axes: $x\in [0,0.5]$.}\label{fig6}
\end{center}
\end{figure}

\begin{figure}[H!]
\begin{center}
\includegraphics[width=11cm]{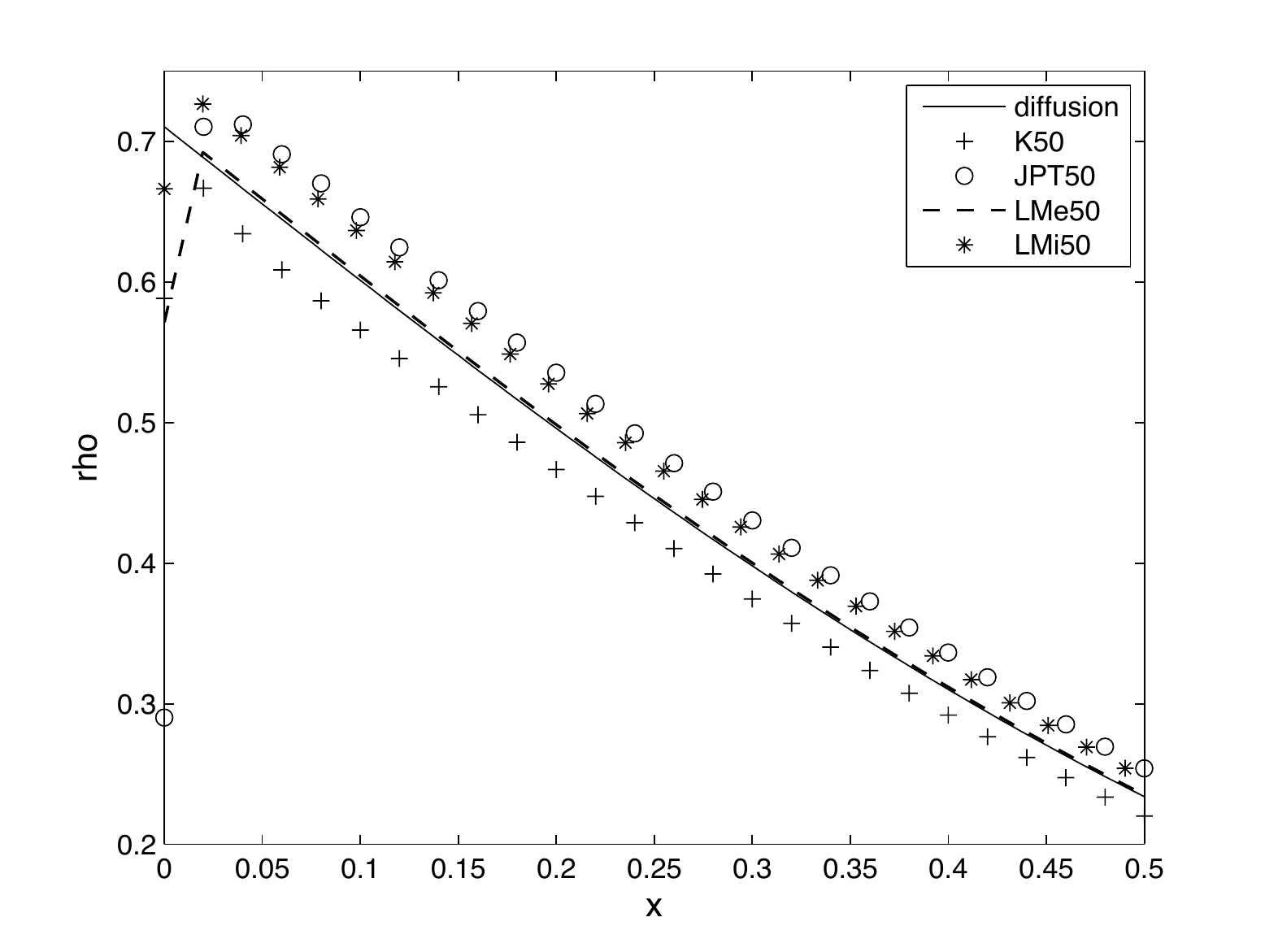}
\caption{Example 5: comparison between the diffusion (200 space gridpoints) and the schemes K, JPT, LMe and LMi (50 space gridpoints each), $\eps=10^{-4}$. A zoom has been made on the axes: $x\in [0,0.5]$.}\label{fig7}
\end{center}
\end{figure}

\begin{figure}[H!]
\begin{center}
\includegraphics[width=11cm]{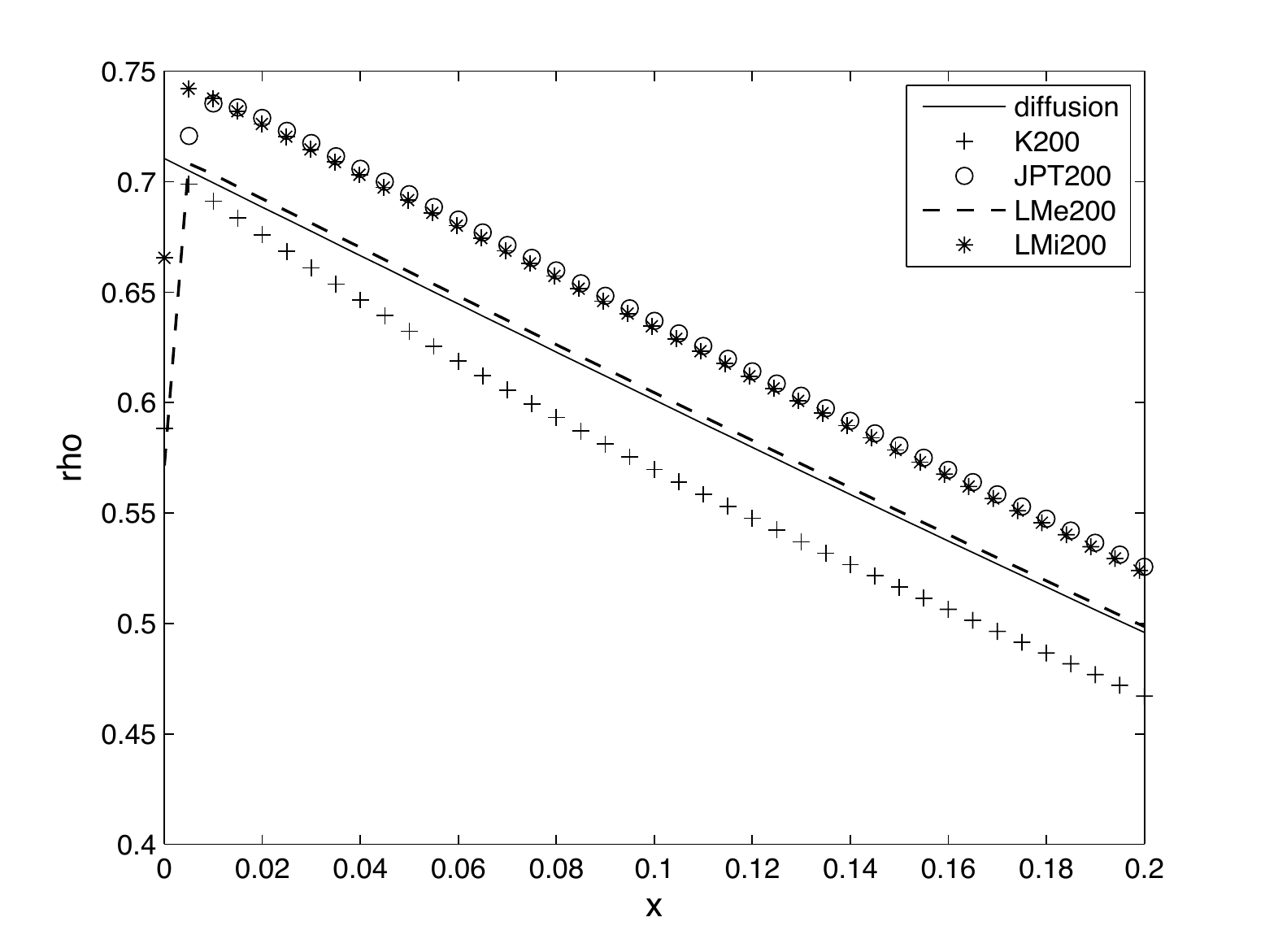}
\caption{Example 5: comparison between the diffusion (200 space gridpoints) and the schemes K, JPT, LMe and LMi (200 space gridpoints each), $\eps=10^{-4}$. A zoom has been made on the axes: $x\in [0,0.2]$.}\label{fig8}
\end{center}
\end{figure}

\begin{figure}[H!]
\begin{center}
\includegraphics[width=11cm]{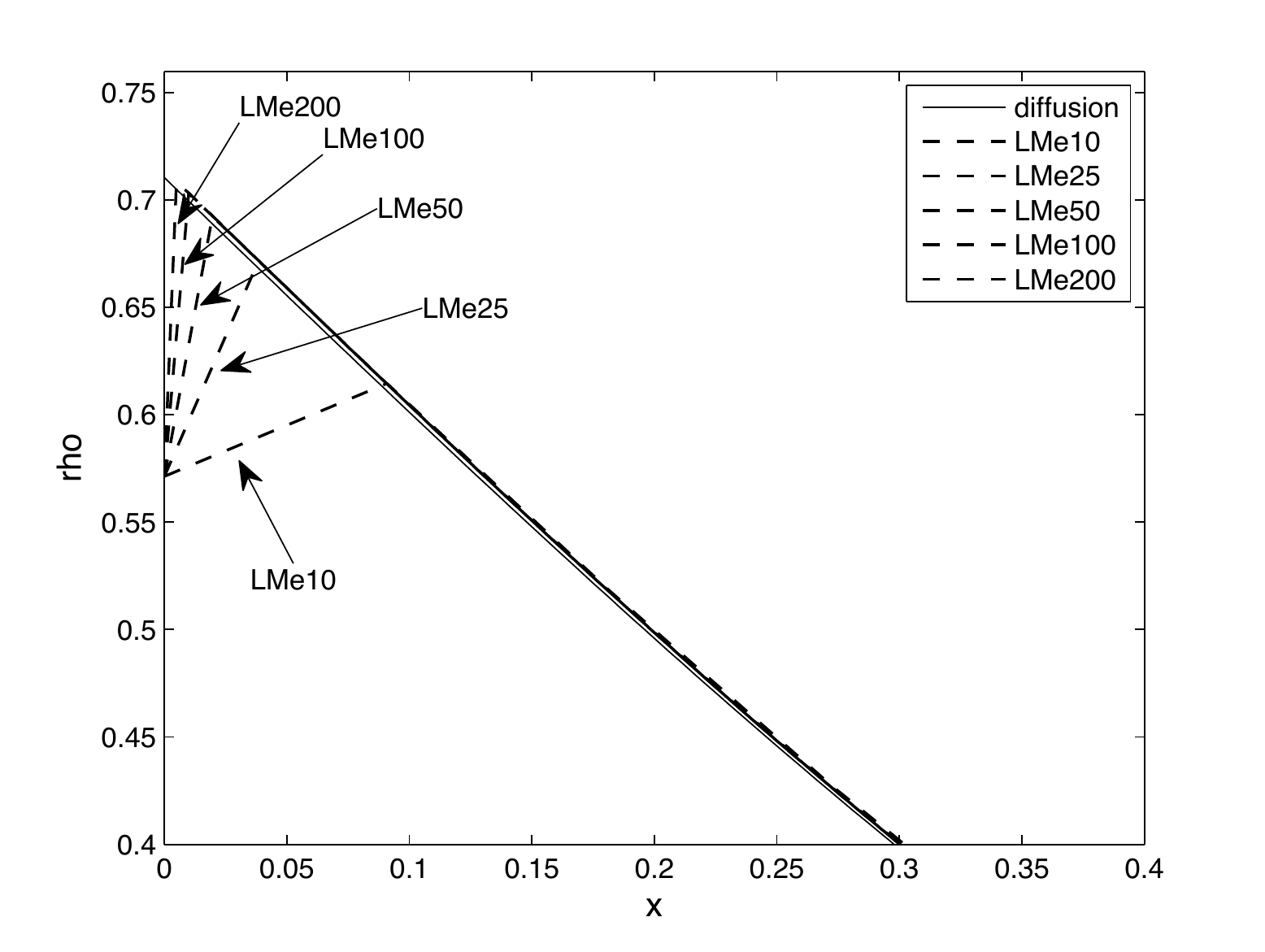}
\caption{Example 5: comparison between the diffusion (200 space gridpoints) and our scheme LMe with 10, 25, 50, 100 and 200 space gridpoints, $\eps=10^{-4}$. A zoom has been made on the axes: $x\in [0,0.4]$.}\label{fig9}
\end{center}
\end{figure}

\begin{figure}[H!]
\begin{center}
\includegraphics[width=11cm]{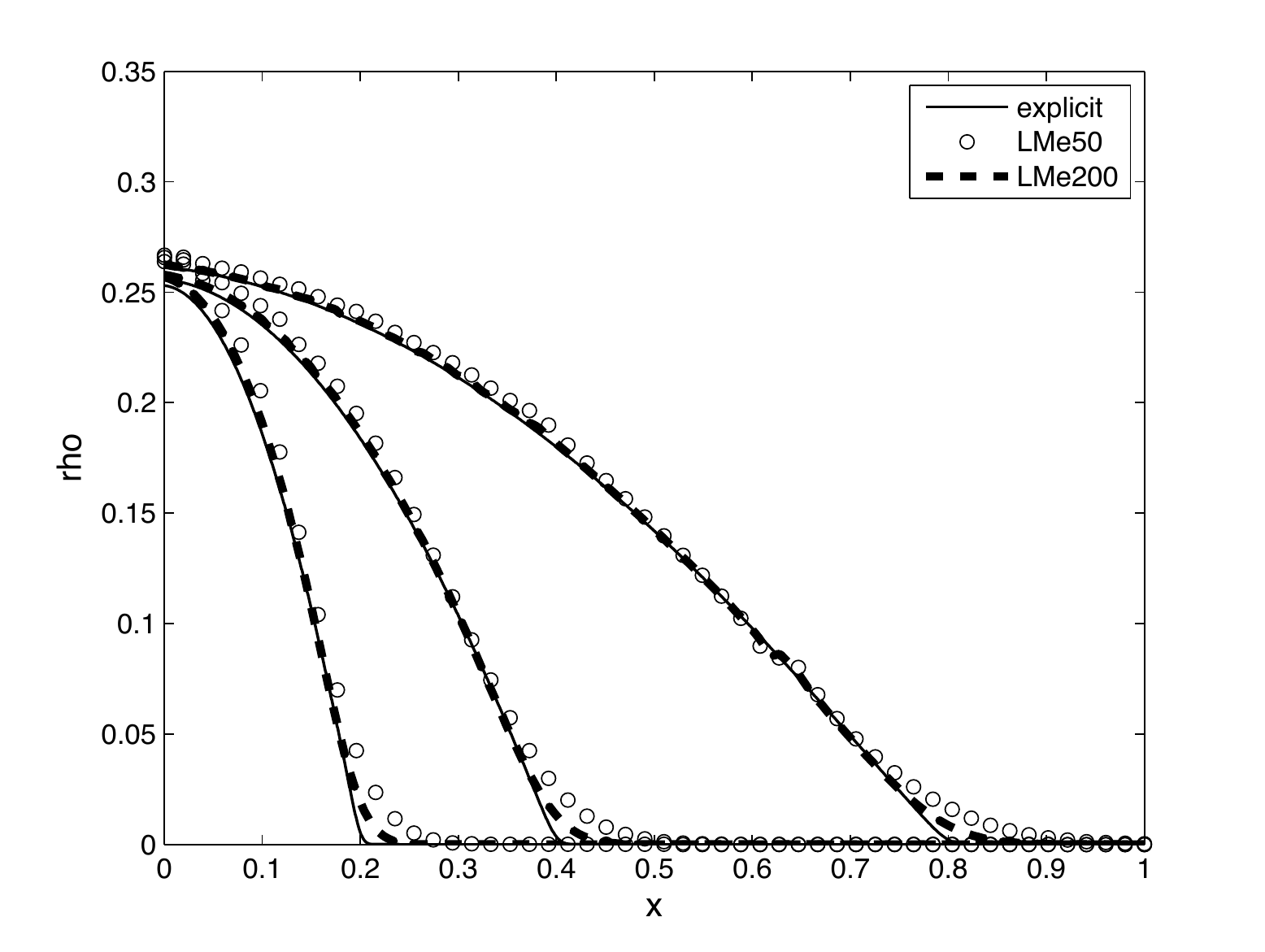}
\caption{Example 6, kinetic regime $\eps=1$: comparison between the explicit scheme (2000 space gridpoints) and our scheme LMe (50 and 200 space gridpoints). Results at times $t=0.2$, $t=0.4$ and $t=0.8$.}\label{fig10}
\end{center}
\end{figure}

\begin{figure}[H!]
\begin{center}
\includegraphics[width=11cm]{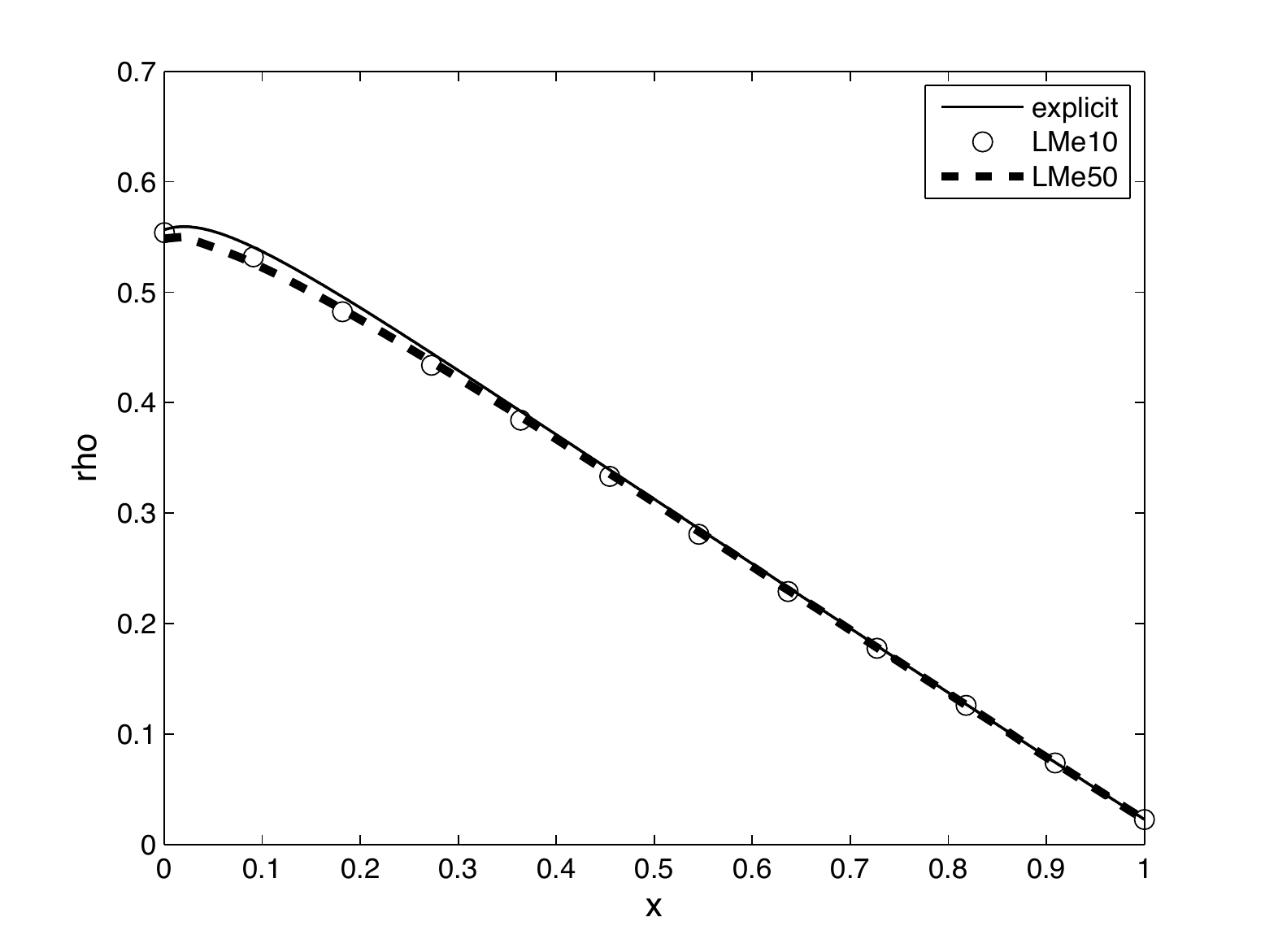}
\caption{Example 6, intermediate regime $\eps=10^{-2}$: comparison between the explicit scheme (10 000 space gridpoints) and our scheme LMe (10 and 50 space gridpoints), at time $t=0.4$.}\label{fig11}
\end{center}
\end{figure}

\begin{figure}[H!]
\begin{center}
\includegraphics[width=11cm]{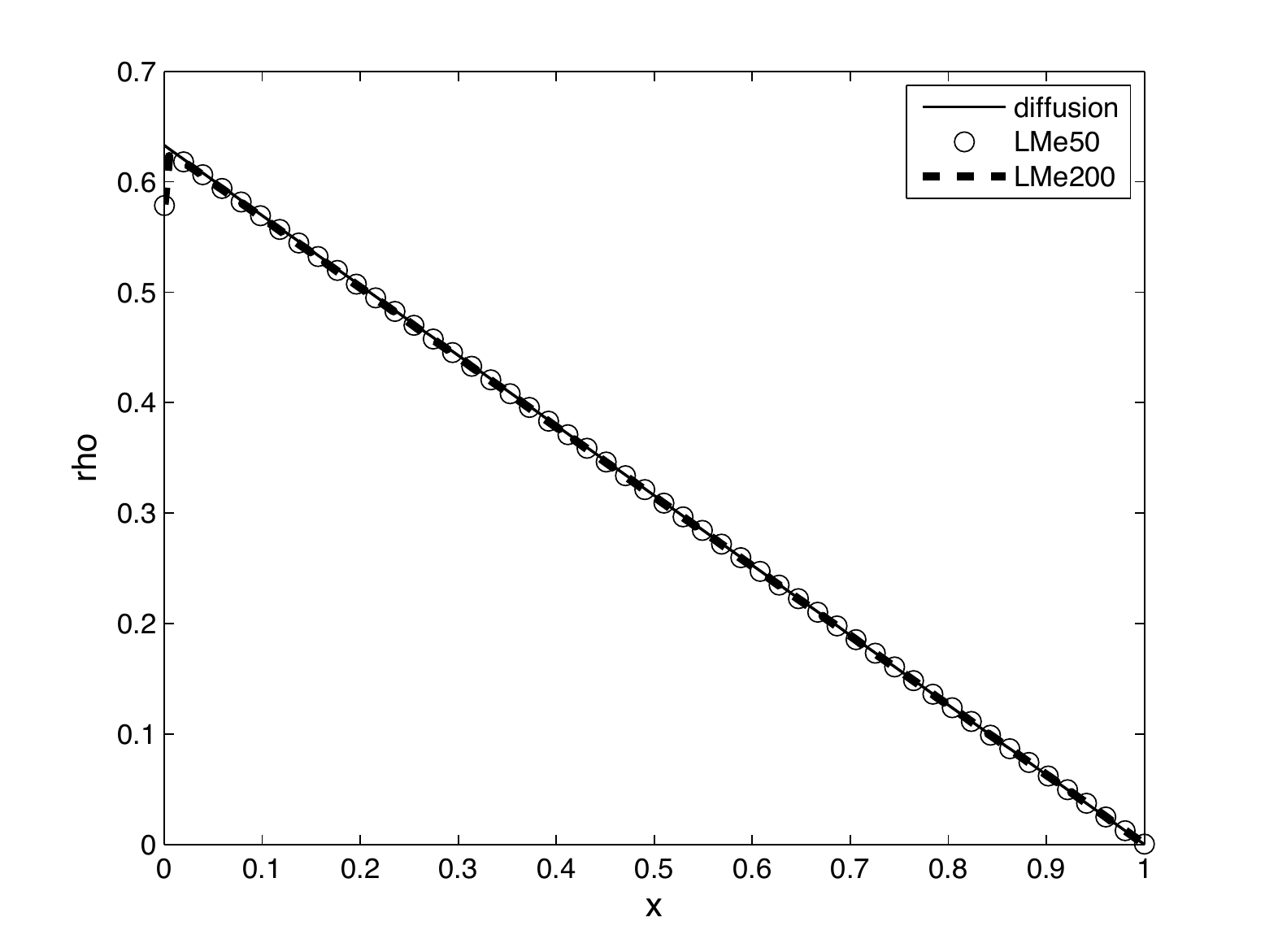}
\caption{Example 6, diffusive regime $\eps=10^{-4}$: comparison between the diffusion and our scheme LMe (50 and 200 space gridpoints), at time $t=0.4$.}\label{fig12}
\end{center}
\end{figure}

\ms
\ni
We plot on Figure \ref{figchandra3} the approximated kernel $K_\sigma^{app}(v)$ given by \fref{polynomesigma2}, compared to the exact kernel $K_\sigma(v)$ computed thanks to \fref{Ksigma}, \fref{chandra2}. We observe again that the kernel generated by our scheme is a very accurate approximation of the exact kernel.

\begin{figure}[H!]
\begin{center}
\includegraphics[width=11cm]{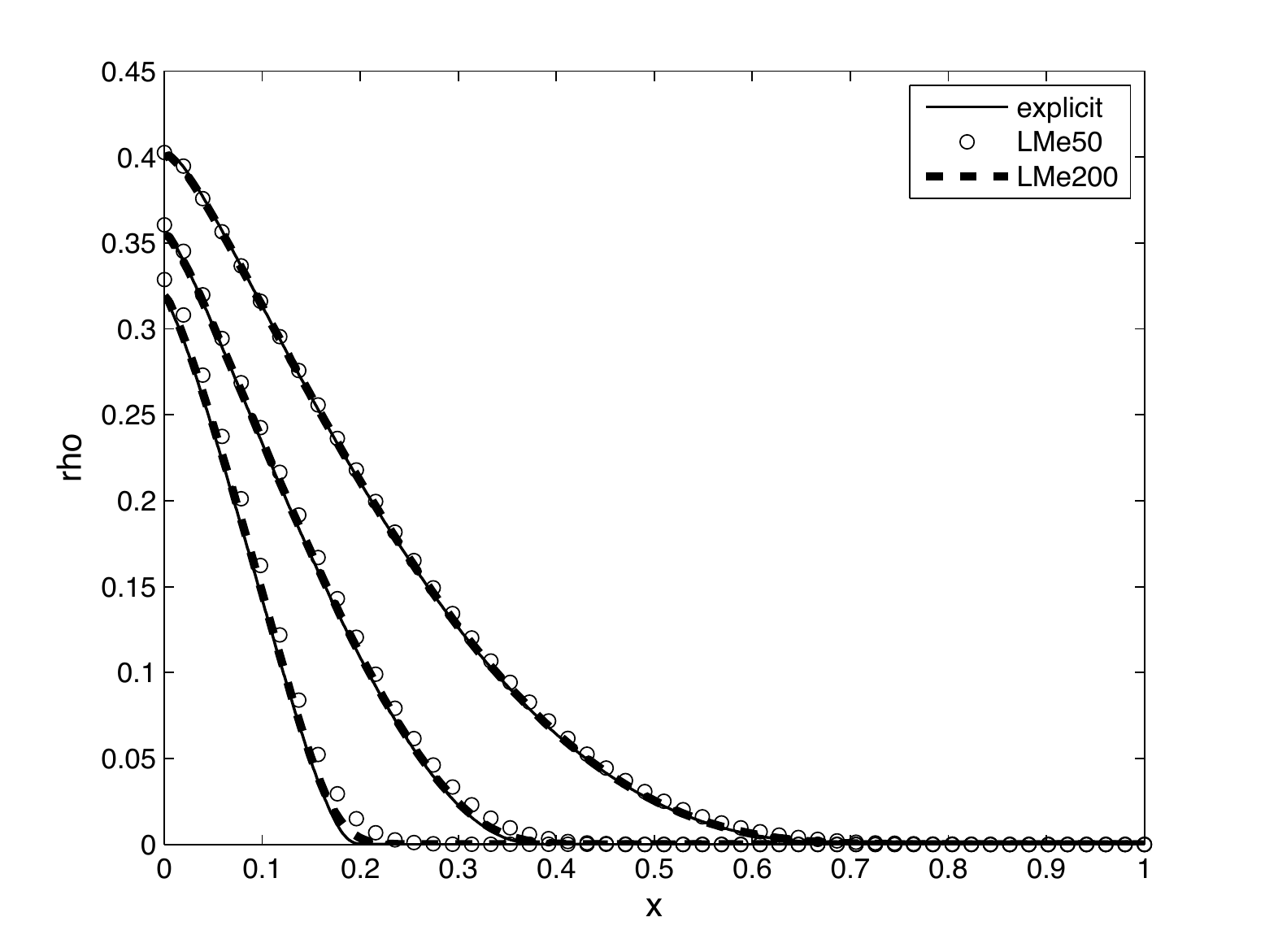}
\caption{Example 7, kinetic regime $\eps=1$: comparison between the explicit scheme (2000 space gridpoints) and our scheme LMe (50 and 200 space gridpoints). Results at times $t=0.2$, $t=0.4$ and $t=0.8$.}\label{fig13}
\end{center}
\end{figure}

\begin{figure}[H!]
\begin{center}
\includegraphics[width=11cm]{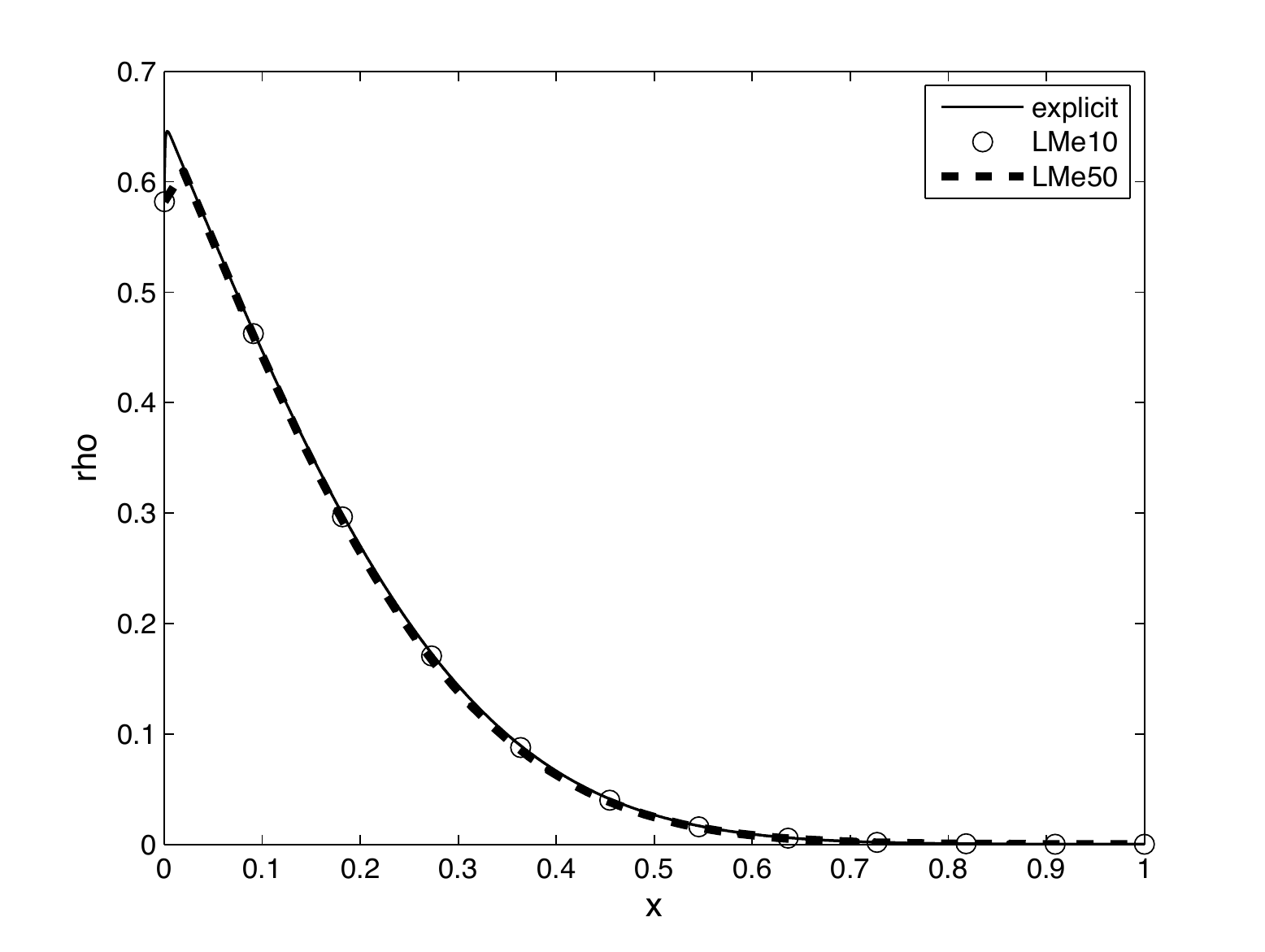}
\caption{Example 7, intermediate regime $\eps=10^{-2}$: comparison between the explicit scheme (10 000 space gridpoints) and our scheme LMe (10 and 50 space gridpoints), at time $t=0.4$.}\label{fig14}
\end{center}
\end{figure}

\begin{figure}[H!]
\begin{center}
\includegraphics[width=11cm]{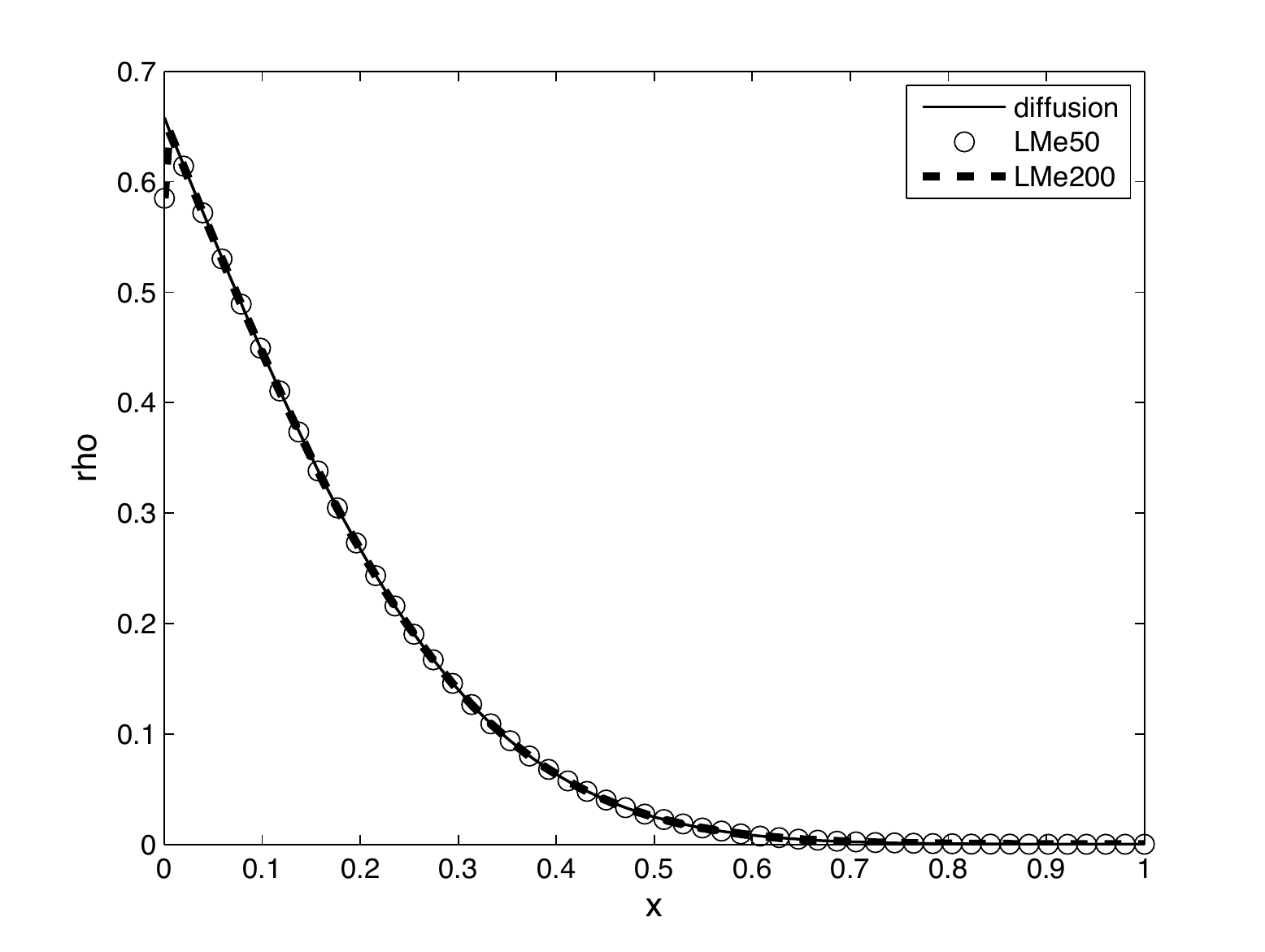}
\caption{Example 7, diffusive regime $\eps=10^{-4}$: comparison between the diffusion and our scheme LMe (50 and 200 space gridpoints), at time $t=0.4$.}\label{fig15}
\end{center}
\end{figure}

We plot the density in the different regimes: $\eps=1$ on Figure \ref{fig13}, $\eps=10^{-2}$ on Figure \ref{fig14} and $\eps=10^{-4}$ on Figure \ref{fig15}. Again, our numerical tests provide a good approximation of the exact result and we do not need any mesh refinement when a boundary layer appears.

\bs
In the two following last examples, we consider cross sections for which formula analogous to \fref{pointfixe} and \fref{chandra2} are not available and could be more complicated to obtain. We only test the intermediate regime $\eps=10^{-2}$, where the reference solution is computed thanks to the highly resolved explicit scheme.

\begin{figure}[H!]
\begin{center}
\includegraphics[width=11cm]{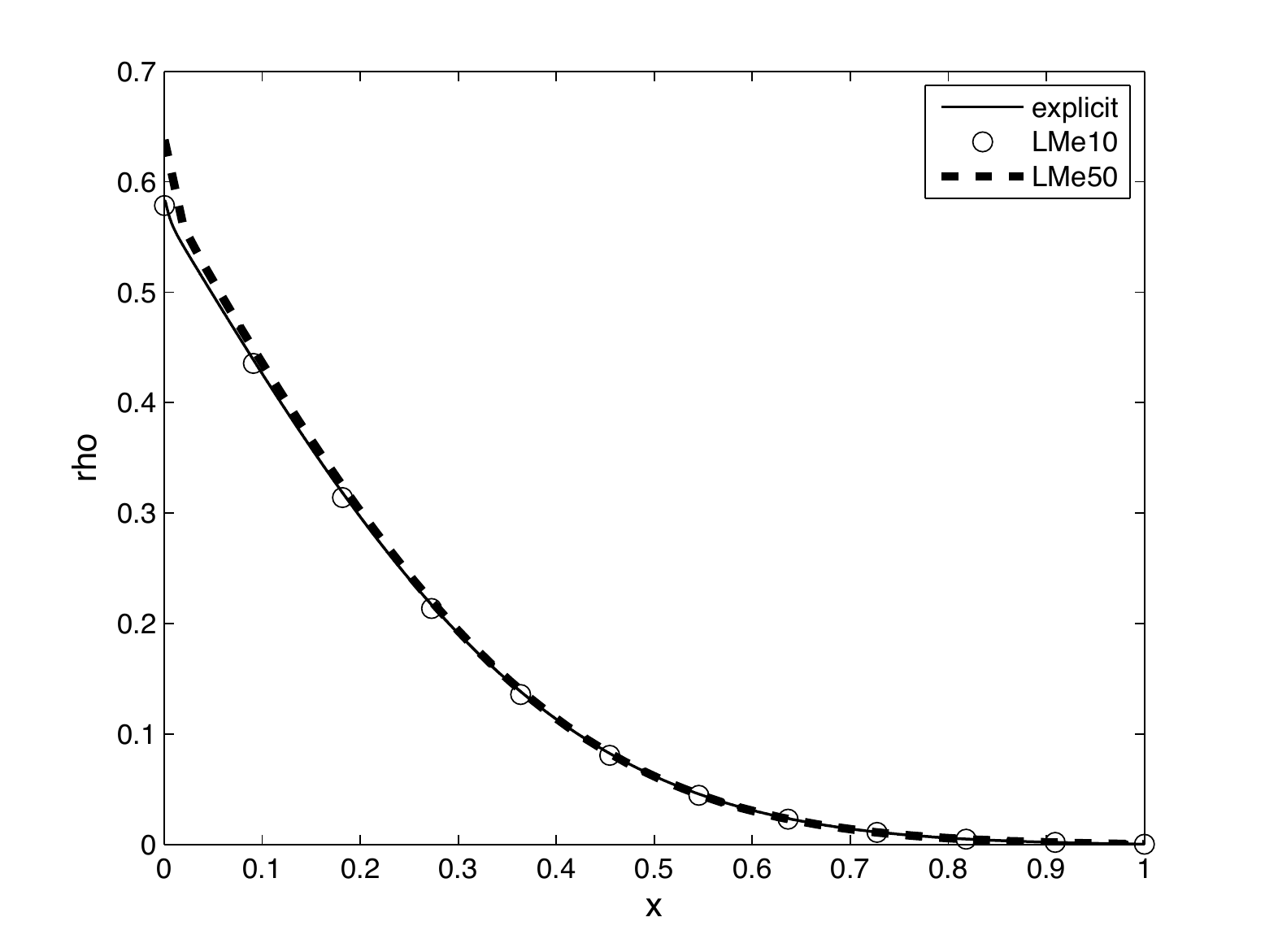}
\caption{Example 8, intermediate regime $\eps=10^{-2}$: comparison between the explicit scheme (10 000 space gridpoints) and our scheme LMe (10 and 50 space gridpoints), at time $t=0.4$.}\label{fig16}
\end{center}
\end{figure}

\bs
\ni
{\em Example 8:} $\sigma(v,w)=|v-w|^{5}$, intermediate regime $\eps=10^{-2}$. The density computed with our scheme is plotted on Figure \ref{fig16} for 10 and 50 space gridpoints and is in good agreement with the reference computation.

\begin{figure}[H!]
\begin{center}
\includegraphics[width=11cm]{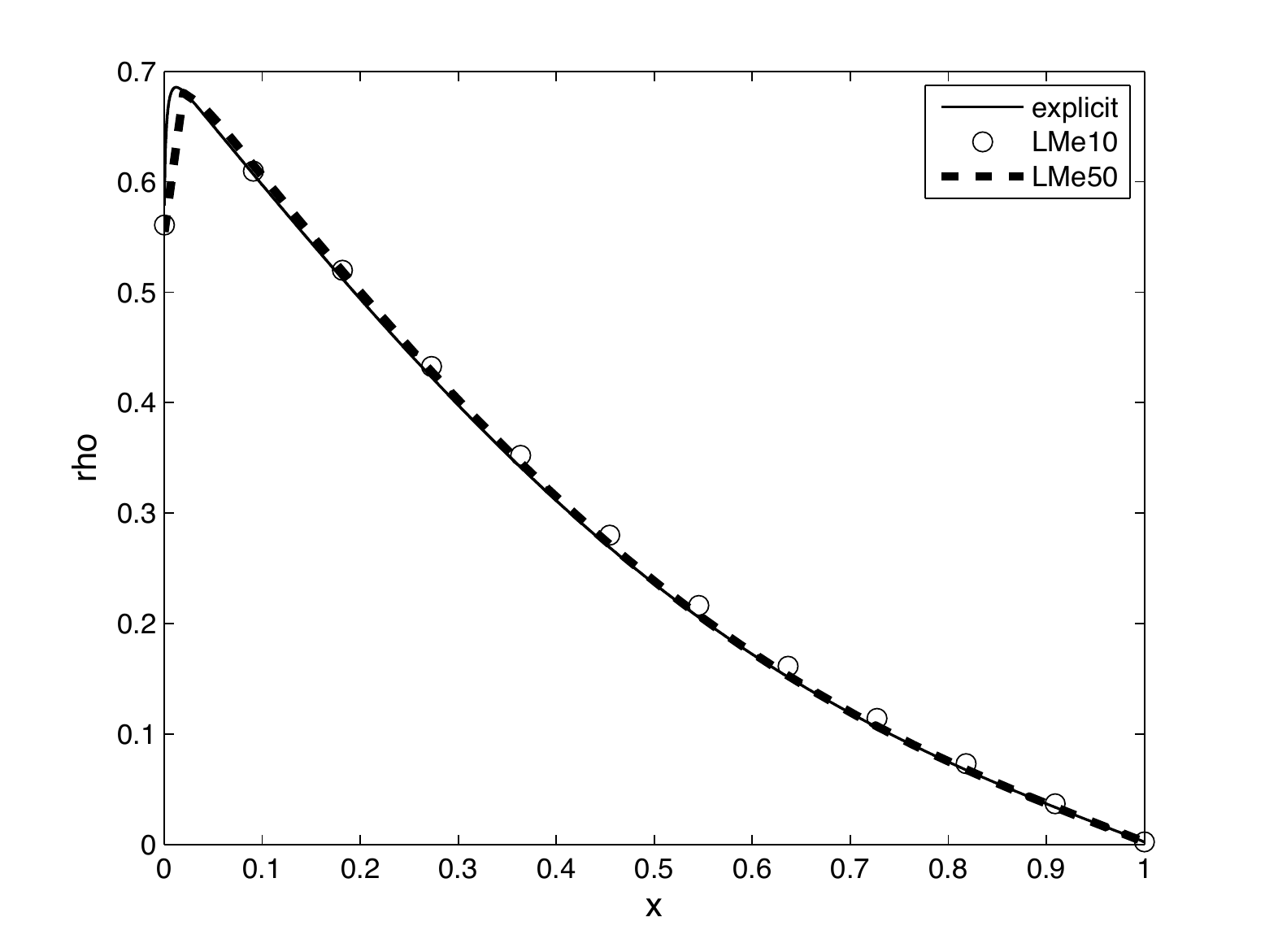}
\caption{Example 9, intermediate regime $\eps=10^{-2}$: comparison between the explicit scheme (10 000 space gridpoints) and our scheme LMe (10 and 50 space gridpoints), at time $t=0.4$.}\label{fig17}
\end{center}
\end{figure}

\bs
\ni
{\em Example 9:} $\sigma(v,w)=|v-w|^{-0.5}$, intermediate regime $\eps=10^{-2}$. The density computed with our scheme is plotted on Figure \ref{fig17}, again for 10 and 50 space gridpoints, and is also in good agreement with the reference computation.

\bs

\section{Conclusion}

In this paper, we have introduced a new strategy to develop the so-called Asymptotic Preserving (AP) scheme for kinetic problems with possible boundary layers in the diffusion limit. This strategy is based on a new micro-macro decomposition. The macroscopic part of the distribution function is not the usual associated Maxwellian but a modified equilibrium which has the same incoming moments. This method is then used to construct a numerical scheme which is AP inside the domain and provides a very good approximation of the exact boundary condition, in both kinetic and diffusion regimes. This scheme only uses the natural inflow boundary condition, and in particular we do not solve the associated Milne problem and we do not need to inject the theoretical limiting boundary value given by the Chandrasekhar type formula (note moreover that such formula is not always available).

We believe that our approach is sufficiently robust to be applied in other situations. For example, we will explore this strategy in the case of hydrodynamic scaling with BGK type collision operators and then with more complex collision operators such as the Landau operator of plasma physics. Moreover, we emphasize that the approach can be extended naturally to multidimensional situations and this task will be achieved is a future work.

\bs
\ni
{\bf Acknowledgements.}
The authors were supported by the Agence Nationale de la Recherche, ANR project QUATRAIN, by the project "D\'efis \'emergeants" funded by the University of Rennes 1 and by the ANR project CBDif. We also acknowledge support from the INRIA team IPSO.


\begin{thebibliography}{9}
\bibitem{BSS} C. Bardos, R. Santos, R. Sentis, Diffusion approximation and computation of the critical size.  Trans. Amer. Math. Soc.  284  (1984),  no. 2, 61--649.
\bibitem{BLM} M. Bennoune, M. Lemou, L. Mieussens,  Uniformly stable numerical schemes for the Boltzmann equation preserving compressible Navier-Stokes asymptotics. J. Comput. Phys. 227 (2008), no. 8, 3781--3803.
\bibitem{BLP} A. Bensoussan, J.-L. Lions, G. C. Papanicolaou,  Boundary layers and homogenization of transport processes.  Publ. Res. Inst. Math. Sci.  15  (1979), no. 1, 53--157.
\bibitem{pareschi-recent} S. Boscarino, L. Pareschi, G. Russo, Implicit-Explicit Runge-Kutta schemes for hyperbolic systems and kinetic equations in the diffusion limit, preprint arXiv:1110.4375.
\bibitem{goudon} J. A. Carrillo, Th. Goudon, P. Lafitte, Simulation of fluid and particles flows: asymptotic preserving schemes for bubbling and flowing regimes. J. Comput. Phys., 227 (2008), no. 16, 7929--7951.
\bibitem{case-zweifel} K. M. Case, P. F.  Zweifel, Linear transport theory. Addison-Wesley Publishing Co., Reading, Mass.-London-Don Mills, Ont. 1967.
\bibitem{chandrasekhar} S. Chandrasekhar, Radiative Transfer. Oxford University Press, 1950.
\bibitem{dautray-lions} R. Dautray, J.-L. Lions,  Analyse math\'ematique et calcul num\'erique pour les sciences et les techniques. Vol. 9. \'Evolution: num\'erique, transport. Masson, Paris, 1988.
\bibitem{degond-masgallic} P. Degond, S. Mas-Gallic, Existence of solutions and diffusion approximation for a model Fokker-Planck equation.  Transport Theory Statist. Phys.  16  (1987),  no. 4-6, 589--636.
\bibitem{dubroca} B. Dubroca, A. Klar, Prise en compte d'un fort d\'es\'equilibre cin\'etique par un mod\`ele aux demi-moments [Half-moment model taking into account strong kinetic non-equilibrium]. C. R. Math. Acad. Sci. Paris  335 (2002), no. 8, 699--704.
\bibitem{jin-filbet} F. Filbet, S. Jin, A class of asymptotic-preserving schemes for kinetic equations and related problems with stiff sources. J. Comput. Phys. 229 (2010), no. 20, 7625--7648.
\bibitem{dubroca2} M. Frank, Martin, B. Dubroca, A. Klar, Partial moment entropy approximation to radiative heat transfer. J. Comput. Phys. 218 (2006), no. 1, 1--18.
\bibitem{golse1} F. Golse, Boundary and interface layers for kinetic problems. Lecture Notes of the 4th Summer school of the GdR SPARCH, 15--20 September 1997, St-Pierre d'Ol\'eron, France.
\bibitem{golse} F. Golse, The Milne problem for the radiative transfer equations (with frequency dependence).  Trans. Amer. Math. Soc.  303  (1987),  no. 1, 125--143.
\bibitem{GJL} F. Golse, S. Jin, C. D. Levermore, The convergence of numerical transfer schemes in diffusive regimes. I. Discrete-ordinate method.  SIAM J. Numer. Anal.   36  (1999),  no. 5, 1333-1369.
\bibitem{gosse-toscani} L. Gosse, G. Toscani,
Asymptotic-preserving \& well-balanced schemes for radiative transfer and the Rosseland approximation.
Numer. Math. 98 (2004), no. 2, 223--250. 
\bibitem{guermond-kanschat} J.-L. Guermond, G. Kanschat, Asymptotic analysis of upwind discontinuous Galerkin approximation of the radiative transport equation in the diffusive limit. SIAM J. Numer. Anal. 48 (2010), no. 1, 53--78.
\bibitem{jin} S. Jin, Efficient asymptotic-preserving (AP) schemes for some multiscale kinetic equations. SIAM J. Sci. Comput. 21 (1999), no. 2, 441--454.
\bibitem{jin-levermore1} S. Jin, D. Levermore,  The discrete-ordinate method in diffusive regimes.  Transport Theory Statist. Phys.  20  (1991),  no. 5-6, 413--439.
\bibitem{jin-levermore2} S. Jin, D. Levermore,  Fully discrete numerical transfer in diffusive regimes.  Transport Theory Statist. Phys.  22  (1993),  no. 6, 739--791.
\bibitem{jin-pareschi} S. Jin, L. Pareschi, Discretization of the multiscale semiconductor Boltzmann equation by diffusive relaxation schemes.  J. Comput. Phys.  161  (2000),  no. 1, 312--330.
\bibitem{JPT} S. Jin, L. Pareschi, G. Toscani, Uniformly accurate diffusive relaxation schemes for multiscale transport equations. SIAM J. Num. Anal., 38 (2000), 913--936.
\bibitem{jin-tang} S. Jin, M. Tang, H. Han, A uniformly second order numerical method for the one-dimensional discrete-ordinate transport equation and its diffusion limit with interface. Netw. Heterog. Media 4 (2009), no. 1, 35--65.
\bibitem{Klar1} A. Klar, An asymptotic-induced scheme for nonstationary transport equations in the diffusion limit. SIAM J. Num.  Anal. 35 (1998), 1073--1094.
\bibitem{Klar3} A. Klar, Asymptotic-induced domain decomposition methods for kinetic and drift diffusion semiconductor equations.  SIAM J. Sci. Comput.  19  (1998),  no. 6, 2032--2050.
\bibitem{klar4} A. Klar, A numerical method for kinetic semiconductor equations in the drift-diffusion limit.  SIAM J. Sci. Comput.  20  (1999),  no. 5, 1696--1712.
\bibitem{Klar2} A. Klar, An asymptotic preserving numerical  scheme for kinetic equations in the low Mach number limit. SIAM J. Num. Anal. 36(1999), 1507--1527.
\bibitem{larsen-keller} E. W. Larsen, J. B.  Keller, Asymptotic solution of neutron transport problems for small mean free paths.  J. Mathematical Phys.  15  (1974), 75--81. 
\bibitem{lafitte} P. Lafitte, G. Samaey, Asymptotic-preserving projective integration schemes for kinetic equations in the diffusion limit, to appear in SIAM J. Sci. Comput.
\bibitem{larsen-morel} E. W. Larsen, J. E. Morel, Asymptotic solutions of numerical transport problems in optically thick, diffusive regimes. II.  J. Comput. Phys.  83  (1989),  no. 1, 212--236.
\bibitem{larsen-morel-miller} E. W. Larsen, J. E. Morel, J. E.; W. F. Miller, Jr, Asymptotic solutions of numerical transport problems in optically thick, diffusive regimes.  J. Comput. Phys.  69  (1987),  no. 2, 283--324.
\bibitem{note} M. Lemou, F. M\'ehats, A boundary matching  micro/macro decomposition for kinetic equations. C. R. Acad. Sci. Paris 349 (2011) 47--484.
\bibitem{LM-AP1} M. Lemou, L. Mieussens,  A new asymptotic preserving scheme based on micro-macro formulation for linear kinetic equations in the diffusion limit. SIAM J. Sci. Comp. 31 (2008), no.1, 334--368.
\bibitem{Lem-note1} M. Lemou,  Relaxed micro-macro schemes for kinetic equations.  C. R. Math. Acad. Sci. Paris   348  (2010),  no. 7-8, 45--460.
\bibitem{liu-mieussens} J.-G. Liu, L. Mieussens, Analysis of an asymptotic preserving scheme for linear kinetic equations in the diffusion limit. SIAM J. Numer. Anal. 48 (2010), no. 4, 1474--1491. 
\bibitem{poupaud} F. Poupaud, Diffusion approximation of the linear semiconductor Boltzmann equation: analysis of boundary layers.  Asymptotic Anal.   4  (1991),  no. 4, 293--317. 
\end{thebibliography}
\end{document}